\documentclass[a4paper]{amsart}

\usepackage[english]{babel}
\usepackage[utf8x]{inputenc}
\usepackage[T1]{fontenc}
\usepackage{subdepth}
\usepackage{cite}

\usepackage{amsthm,amsmath,amssymb}
\usepackage{tikz-cd}
\usepackage{amsfonts}
\usepackage{enumitem}
\usepackage{graphicx}
\usepackage[colorlinks=true, allcolors=blue]{hyperref}
\usepackage{cleveref}
\usepackage{mathdots}
\usepackage{mathtools}
\usepackage{csquotes}
\usepackage{xfrac}
\usepackage{ulem}
\usepackage{amsaddr}

\tikzset{no body/.style={/tikz/dash pattern=on 0 off 1mm}}

\usetikzlibrary{decorations.markings,arrows,patterns}

\theoremstyle{plain}
\newtheorem{thm}{Theorem}[section]
\newtheorem{lem}[thm]{Lemma}
\newtheorem{prop}[thm]{Proposition}
\newtheorem{cor}[thm]{Corollary}
\newtheorem*{thm*}{Theorem}
\numberwithin{equation}{section}
\newtheorem{thmx}{Theorem}

\newtheorem{corx}[thmx]{Corollary}
\theoremstyle{definition}
\newtheorem{example}[thm]{Example}

\newtheorem{defn}[thm]{Definition}

\theoremstyle{remark}

\newtheorem{rmk}[thm]{Remark}

\newcommand*{\suchthat}{\;\ifnum\currentgrouptype=16 \middle\fi|\;}
\newcommand{\Z}{\mathbb{Z}}
\newcommand{\I}{\mathfrak{I}}
\newcommand{\T}{\mathfrak{T}}
\DeclareMathOperator{\KK}{KK}

\newcommand{\acts}{\mathbin{\curvearrowright}}
\newcommand{\rightacts}{\mathbin{\curvearrowleft}}

\newcommand{\A}{\mathcal{A}}
\newcommand{\B}{\mathcal{B}}

\renewcommand{\H}{\mathrm{H}}
\newcommand{\Ho}{\mathrm{Ho}}
\newcommand{\K}{\mathrm{K}}
\newcommand{\Ktop}{\mathrm{K}^{\mathrm{top}}}
\newcommand{\Kred}{\mathrm{K}^{\mathrm{red}}}
\DeclareMathOperator{\Ab}{\mathbf{Ab}}

\newcommand{\lMod}[1]{#1\text{-}\mathbf{Mod}}

\newcommand{\cat}[1]{\mathbf{#1}}
\newcommand{\bb}[1]{\mathbb{#1}}
\newcommand{\fk}[1]{\mathfrak{#1}}
\DeclareMathOperator{\ABC}{ABC}

\DeclareMathOperator{\EC}{EC}
\DeclareMathOperator{\App}{App}

\newcommand{\cal}[1]{\mathcal{#1}}
\newcommand{\bcdot}{\boldsymbol{\cdot}}
\renewcommand{\restriction}{\mathord{\upharpoonright}}
\DeclareMathOperator{\im}{im}
\DeclareMathOperator{\id}{id}
\DeclareMathOperator{\orb}{orb}
\DeclareMathOperator{\dom}{dom}

\DeclareMathOperator{\Ind}{Ind}
\DeclareMathOperator{\Res}{Res}

\DeclareMathOperator{\Cone}{Cone}

\newcommand{\downthenrightarrow}{\tikz{\draw[->] (0,0.2) |- (0.25,0);}}
\newcommand{\rightthendownarrow}{\tikz{\draw[->] (0,0.25) -| (0.25,0);}}
\newcommand{\defeq}{:=}
\newcommand{\orho}{\overline{\rho}}
\newcommand{\cs}{\mathrm{C}^\ast}

\makeatletter
\tikzcdset{
  open/.code     = {\tikzcdset{hook, circled};},
  closed/.code   = {\tikzcdset{hook, slashed};},
  open'/.code    = {\tikzcdset{hook', circled};},
  closed'/.code  = {\tikzcdset{hook', slashed};},
  circled/.code  = {\tikzcdset{markwith = {\draw[black, fill = white] (0,0) circle (.375ex);}};},
  circledred/.code  = {\tikzcdset{markwith = {\draw[red, fill = white] (0,0) circle (.375ex);}};},
  circledblue/.code  = {\tikzcdset{markwith = {\draw[blue, fill = white] (0,0) circle (.375ex);}};},
  slashed/.code  = {\tikzcdset{markwith = {\draw[-] (-.4ex,-.4ex) -- (.4ex,.4ex);}};},
  markwith/.code ={
    \pgfutil@ifundefined%
    {tikz@library@decorations.markings@loaded}%
    {\pgfutil@packageerror{tikz-cd}{You need to say %
      \string\usetikzlibrary{decorations.markings} to use arrows with markings}{}}{}%
    \pgfkeysalso{/tikz/postaction = {
      /tikz/decorate,
      /tikz/decoration={markings, mark = at position 0.5 with {#1}}}
    }
  },
}
\makeatother

\newcommand{\circledrightarrow}{\mathrel{\begin{tikzpicture}[
    decoration={markings,mark=at position 0.5 with {{\draw[black, fill = white] (0,0) circle (.25ex);}}},
    ]
    \draw[postaction={decorate}, ->] (0,0) -- (0.34,0);
\end{tikzpicture}}}

\newcommand{\xcircledrightarrow}[2][->]{
\tikz[baseline=-0.6ex,
    decoration={markings,mark=at position 0.5 with {{\draw[black, fill = white] (0,0) circle (.25ex);}}},
    ]{
\node[anchor=south,font=\scriptsize, inner ysep=2.8pt,outer xsep=2.5pt](x){$#2$};
\draw[shorten <=3.4pt,shorten >=3.4pt,postaction=decorate,#1](x.south west)--(x.south east);
}
}

\newcommand\extrafootertext[1]{%
    \bgroup
    \renewcommand\thefootnote{\fnsymbol{footnote}}%
    \renewcommand\thempfootnote{\fnsymbol{mpfootnote}}%
    \footnotetext[0]{#1}%
    \egroup
}

\title[Isomorphisms in K-theory from isomorphisms in groupoid homology]{Isomorphisms in K-theory from isomorphisms in groupoid homology theories}
\author{Alistair Miller}
\address{}

\begin{document}
\maketitle

\begin{abstract}

We prove that for torsion-free amenable ample groupoids, an isomorphism in groupoid homology induced by an étale correspondence yields an isomorphism in the K-theory of the associated $\mathrm{C}^\ast$-algebras. We apply this to extend X. Li's K-theory formula for left regular inverse semigroup $\mathrm{C}^\ast$-algebras. These results are obtained by developing the functoriality of the ABC spectral sequence.

\end{abstract}

\extrafootertext{This work contains part of the author's PhD thesis, which was supported by the Engineering and Physical Sciences Research Council (EPSRC) through a doctoral studentship. The author has also been supported by the European Research Council (ERC) under the European Union's Horizon 2020 research and innovation programme (grant agreement No. 817597) and by the Independent Research Fund Denmark through the Grant 1054-00094B.}

\section{Introduction}

Étale groupoids and their $\cs$-algebras are a central topic in operator algebras, arising naturally from noncommutative geometry and topological dynamics. Twisted étale groupoids also arise from within the field of $\cs$-algebras via reconstruction \cite{Renault08} and as models for every Elliott-classifiable $\cs$-algebra \cite{Li20}. The homological and K-theoretic invariants associated to étale groupoids contain crucial information. Indeed, the power of K-theory to remember the associated $\cs$-algebra is demonstrated by the celebrated classification programme \cite{Kirchberg95, Phillips00, GLN20, TWW17, Winter18}, and the K-theory furthermore often retains information prior to the $\cs$-algebra, for example classifying Cantor minimal systems up to orbit equivalence \cite{GPS95}. While the computation of the operator K-theory $\K_*(C^*_r(G))$ associated to an étale groupoid $G$ can prove quite challenging, it is often possible to make progress by relating the operator K-theory to more computable homological invariants of the groupoid such as groupoid homology and topological K-theory. Our main contribution is to show that an isomorphism in groupoid homology can be converted to an isomorphism in K-theory when it is induced by an étale correspondence.

A connection between the homology of an ample groupoid and the operator K-theory was uncovered by Matui in examples coming from topological dynamics \cite{Matui12, Matui15, Matui16}. This led him to conjecture that for an ample groupoid $G$ in a certain class, the K-theory groups decompose as a direct sum of groupoid homology groups of the same parity:
\[ \K_*(C^*_r(G)) \cong \bigoplus_{i \in \bb N} \H_{*+2i}(G) \]
Known as Matui's HK conjecture, this has been verified for a wide variety of examples \cite{FKPS19, BDGW23, ProYam23}, but has counterexamples \cite{Scarparo20, Deeley23} which ultimately point to a more subtle relationship between $\H_*(G)$ and $\K_*(C^*_r(G))$.

The topological K-theory $\Ktop_*(G)$ of a second countable Hausdorff étale groupoid $G$ is an alternative K-theoretic invariant to the operator K-theory. The Baum--Connes conjecture for $G$ \cite{BCH94, Tu99b} asserts the bijectivity of a homomorphism 
\[ \mu_G \colon \Ktop_*(G) \to \K_*(C^*_r(G)) \]
known as the Baum--Connes assembly map. This is known to hold for groupoids with the Haagerup property (in particular amenable groupoids) \cite{Tu99b}, and while there are counterexamples \cite{HLS02}, the conjecture is open for exact groupoids. The nature of topological K-theory makes it more susceptible to homological methods. Specialising to groups, there is a spectral sequence which converges to the topological K-theory whose second sheet is given by Bredon homology \cite{Mislin03}. An instance of Meyer's ABC spectral sequence \cite{Meyer08} considered by Proietti and Yamashita for ample groupoids with torsion-free isotropy converges to the topological K-theory, with groupoid homology on the second sheet \cite{ProYam22}. Due to work of B\"onicke and Proietti \cite{BoePro24}, there is a more general spectral sequence converging to the topological K-theory of an arbitrary étale groupoid:
\begin{equation*}\label{intro spec seq}
E^2_{p,q} \Rightarrow \Ktop_{p+q}(G)
\end{equation*}
The second sheet of this spectral sequence is occupied by a homology theory for the groupoid analogous to Bredon homology for groups and groupoid homology for torsion-free ample groupoids. Outside of the torsion-free ample setting we should expect a homology theory which differs from groupoid homology. This sheds some light on the HK conjecture; the presence of torsion in the isotropy makes groupoid homology the wrong invariant to compare to K-theory (see \cite{Scarparo20}), and while the spectral sequence degenerates to a direct sum when the homological dimension is low, this should not be expected in higher dimensions (see \cite{Deeley23}).

To get the most out of the above spectral sequences linking these homology theories to operator/topological K-theory, we study their functoriality with respect to the groupoid. We work with a class of morphisms of étale groupoids called étale correspondences, which are an analogue for étale groupoids of $\cs$-correspondences, modelling the behaviour of both $*$-homomorphisms and Morita equivalences. They encompass many kinds of morphisms of étale groupoids, including open embeddings, Morita equivalences, actors\footnote{also known as Zakrzewski or algebraic morphisms} and étale homomorphisms. A proper étale correspondence $\Omega \colon G \to H$ induces a map in K-theory through the associated proper $\cs$-correspondence \cite{AKM22}, and for ample groupoids the author constructed an induced map $\H_*(\Omega) \colon \H_*(G) \to \H_*(H)$ in homology \cite{Miller24}.

While the situation is more complex outside of the torsion-free ample setting, the spirit of this paper is captured by the following result:

\begin{thmx}[see \Cref{isomorphism of homology into k theory}]\label{isomorphism theorem intro}
Let $G$ and $H$ be second countable Hausdorff ample groupoids with torsion-free isotropy groups and which satisfy the Baum--Connes conjecture. Suppose $\Omega \colon G \to H$ is a proper second countable Hausdorff étale correspondence which induces an isomorphism
\[\H_{*}(\Omega) \colon \H_*(G) \to \H_*(H)\]
in homology. Then $\Omega$ induces an isomorphism $\K_*(C^*_r(G)) \cong \K_*(C^*_r(H))$ in K-theory. 
\end{thmx}
If furthermore $\Omega \colon G \to H$ induces a $\cs$-correspondence $C^*_r(\Omega) \colon C^*_r(G) \to C^*_r(H)$ at the reduced level (see \Cref{functoriality of topological K theory}), the induced isomorphism is given by 
\[\K_*(C^*_r(\Omega)) \colon \K_*(C^*_r(G)) \to \K_*(C^*_r(H)).\]

Groupoid homology is often more tractable than operator K-theory, and it may be possible to demonstrate an isomorphism in homology directly where there is no clear direct method for K-theory. For example, for any inverse semigroup $S$ with universal groupoid $G_S$, there is a proper étale correspondence 
\[ \Omega_S \colon S^\times \to G_S, \]
viewing the set $S^\times$ of non-zero elements of $S$ as a discrete groupoid. This induces an isomorphism in homology \cite[Example 3.10]{Miller24}. Applying \Cref{isomorphism theorem intro} to $\Omega_S$ enables us to reduce the computation of the K-theory of the reduced $\cs$-algebra $C^*_\lambda(S)$ to computations for its stabiliser subgroups at non-zero idempotents: 
\begin{corx}[see \Cref{inverse semigroup corollary}]\label{intro inverse semigroup cor}
Let $S$ be a countable inverse semigroup with idempotent semilattice $E$ whose universal groupoid is Hausdorff and satisfies the Baum--Connes conjecture. Suppose further that the stabiliser subgroup $S_e = \{ s \in S \mid s^* s = ss^* = e \}$ of each non-zero idempotent $e \in E^\times$ is torsion-free. Then
\[ \K_*(C^*_\lambda(S)) \cong \bigoplus_{\orb(e) \in S \backslash E^\times} \K_*(C^*_\lambda(S_e)). \]
\end{corx}

Aside from our torsion-freeness assumption, this generalises \cite[Theorem 1.1]{Li21a} to inverse semigroups that need not admit an idempotent pure partial homomorphism to a group. Let us explain the key differences in our work to previous approaches \cite{CEL13, CEL15, Li21a} that allow us to overcome the difficulties in extending to this more general setting. When $S$ admits an idempotent pure partial homomorphism to a group $\Gamma$, the dynamics of $S$ are governed by $\Gamma$ in such a way that the desired map $\K_*(C^*_r(S^\times)) \to \K_*(C^*_\lambda(S))$ can be constructed from a $\Gamma$-equivariant Kasparov cycle $\alpha \in \KK^\Gamma$ as $\K_*(\Gamma \ltimes_r \alpha)$. Here the previous approaches make vital use of the going-down principle \cite{CEOO04}, which gives a sufficient condition on an equivariant Kasparov cycle to induce an isomorphism in the K-theory of the crossed products. However, in general there need be no group to govern the dynamics of $S$, so we turn to groupoid techniques. Although a going-down principle for ample groupoids was developed in \cite{Boenicke20}, there is no obvious way to construct the desired Kasparov cycle equivariantly with respect to a single groupoid. Instead, we model the 
desired map in K-theory with the proper étale correspondence $\Omega_S \colon S^\times \to G_S$ of groupoids. The sufficient condition given in \Cref{isomorphism theorem intro} to obtain an isomorphism in K-theory may be compared to the going-down principle, but rather than having to work equivariantly with respect to a single group or groupoid we can work with morphisms on the groupoid level.

\Cref{intro inverse semigroup cor} yields a K-theory formula for the left regular algebra $C^*_\lambda(P)$ of many left cancellative monoids $P$ which satisfy the independence condition but which need not embed in a group (see \Cref{monoid example}). It also recovers the K-theory of the Toeplitz algebra of a general singly aligned higher rank graph (\Cref{HRG example}).

The results in this paper are centred around the ABC spectral sequence \cite[Theorem 5.1]{Meyer08}, which is constructed in the Meyer--Nest framework of homological algebra in triangulated categories via homological ideals \cite{MeyNes06, MeyNes10, Meyer08}. The main example we focus on is in the setting of a second countable Hausdorff étale groupoid $G$ with a specified countable family $\cal F$ of proper open subgroupoids, meant to resemble the family of finite subgroups of a countable group, and a separable $G$-$\cs$-algebra $A$. There is a homological ideal $\I_{\cal F}$ in the (triangulated) equivariant Kasparov category $\KK^G$ consisting of the morphisms which vanish under restriction $\Res^H_G \colon \KK^G \to KK^H$ to every member $H$ of $\cal F$. With respect to the homological ideal $\I_{\cal F}$ we can form the \textit{localisation} $\bb L^{\cal F} \K_*(G \ltimes_r -)$ and the \textit{derived functors} $\bb L^{\cal F}_n \K_*(G \ltimes_r -)$ for $n \geq 0$ of the stable homological functor $\K_*(G \ltimes_r -)$ on $\KK^G$. The localisation is an approximation to $\K_*(G \ltimes_r -)$ which vanishes on the objects in $\I_{\cal F}$, and comes with a natural transformation $\bb L^{\cal F} \K_*(G \ltimes_r -) \Rightarrow \K_*(G \ltimes_r -)$ to the original functor. For suitable $\cal F$ the localisation may be identified with the topological K-theory, through which the natural transformation recovers the Baum--Connes assembly map \cite[Theorem C]{BoePro24}. The derived functors are analogous to classical derived functors on abelian categories; when $G$ is ample and $\cal F = \{G^0\}$ they model the groupoid homology groups $\H_n(G;\K_*(-))$ \cite{ProYam22}. The ABC spectral sequence approximates the localisation of a stable homological functor $F$ on a triangulated category $\T$ at an object $A$ with respect to a suitable homological ideal $\I \vartriangleleft \T$ with its derived functors:
\begin{equation*}
E^2_{p,q} = \bb L^\I_p F_q(A) \Rightarrow \bb L^\I F_{p+q} (A)
\end{equation*}
We refer to such quadruples $(\T,\I,F,A)$ as \textit{ABC tuples} and introduce \textit{ABC morphisms} to place them in a category. We construct ABC morphisms from étale correspondences $\Omega \colon G \to H$ of étale groupoids which are compatible with respect to families of subgroupoids on either side. The two main components of this construction are the \textit{induction functor} $\Ind_\Omega \colon \KK^H \to \KK^G$ and the \textit{induction natural transformation} $\alpha_\Omega \colon \K_*(G \ltimes \Ind_\Omega -) \Rightarrow \K_*(H \ltimes -)$ from \cite{Miller23a}. From an ABC morphism $\fk m$ we construct a morphism $\bb L(\fk m)$ of the associated localisations, morphisms $\bb L_n(\fk m)$ of the derived functors and ultimately a morphism $\ABC(\fk m)$ of the associated spectral sequences:
\begin{thmx}[see Theorem \ref{ABCthmintext}]\label{ABC thm}
An ABC morphism $\fk m \colon \fk M \to \fk M'$ functorially induces a morphism of ABC spectral sequences $\ABC(\fk m) \colon \ABC(\fk M) \to \ABC(\fk M')$. Moreover,
\begin{enumerate}[label=(\roman*)]
\item the map on the second sheet is given by the derived functor maps $\bb L_n(\fk m)$, 
\item the map on the limit sheet agrees with the localisation map $\bb L(\fk m)$.
\end{enumerate}
\end{thmx}
An ABC morphism $\fk m$ for which the morphism $\bb L_n(\fk m)$ of derived functors is an isomorphism for each $n$ therefore induces an isomorphism $\bb L(\fk m)$ of localisations \cite[Comparison Theorem 5.2.12]{Weibel94}. Since we wish to apply this principle to the topological K-theory of étale groupoids, we initiate a systematic study into ABC tuples we can associate to a second countable Hausdorff étale groupoid $G$ whose localisation models the topological K-theory. While suitable families of proper open subgroupoids have been constructed in order to obtain the topological K-theory in \cite{BoePro24}, there is no guarantee that a given étale correspondence is compatible with families of this form on either side. We seek to remedy this through the introduction of a practical sufficient condition for a family to be suitable. A family $\cal F$ of proper open subgroupoids of $G$ satisfies \textit{condition (P)} if for each $x \in G^0$ and each finite subgroup $\Gamma \leq G^x_x$ of the isotropy group at $x$, there is a member $H \in \cal F$ containing $\Gamma$. Immediately it is clear that if $G$ has torsion-free isotropy, then the trivial family $\{G^0\}$ satisfies condition (P). In \Cref{inverse semigroup P example} we introduce an explicit way to construct families satisfying (P) for any étale groupoid.

\begin{thmx}[see Propositions \ref{condition P gives topological K theory} and \ref{checkable P}]\label{intro theorem for obtaining topological K theory}
Let $G$ be a second countable Hausdorff étale groupoid and suppose a countable family $\cal F$ of proper open subgroupoids of $G$ satisfies condition (P). Then for $A \in \KK^G$ there is a natural isomorphism
\[ \bb L^{\cal F} \K_*(G \ltimes_r A) \cong \Ktop_*(G;A).   \]
\end{thmx}

This can be pushed further if we relax the notion of an open subgroupoid to a groupoid which need only be openly embedded up to Morita equivalence, which we term an \textit{open Morita embedding}. The construction of homological ideals giving ABC tuples generalises to families of open Morita embeddings, as does condition (P) and subsequently \Cref{intro theorem for obtaining topological K theory}. We describe constructions of families of open Morita embeddings satisfying condition (P) compatible with a given étale correspondence, yielding:

\begin{thmx}[see \Cref{functoriality of topological K theory}]\label{intro theorem on localisation functoriality}
Let $G$ and $H$ be second countable Hausdorff étale groupoids with coefficient $\cs$-algebras $A \in \KK^G$ and $B \in \KK^H$. Given a second countable Hausdorff étale correspondence $\Omega \colon G \to H$ and a morphism $f \in \KK^G(A, \Ind_\Omega B)$ there is an induced map
\[ \Ktop_*(\Omega;f) \colon \Ktop_*(G;A) \to \Ktop_*(H;B)\]
in topological K-theory. Moreover, $\Ktop_*(\Omega;f)$ may be identified with the map $\bb L(\fk m)$ of localisations induced by an ABC morphism $\fk m$.
\end{thmx}

Without coefficients, this specialises to a map $\Ktop_*(\Omega) \colon \Ktop_*(G) \to \Ktop_*(H)$ induced by a \textit{proper} étale correspondence $\Omega \colon G \to H$. Combining Theorems \ref{ABC thm} and \ref{intro theorem on localisation functoriality}, we have shown that to prove $\Ktop_*(\Omega;f) \colon \Ktop_*(G;A) \to \Ktop_*(H;B)$ to be an isomorphism, it suffices to construct an ABC morphism $\fk m$ with $\bb L(\fk m) = \Ktop_*(\Omega;f)$ such that $\bb L_n(\fk m)$ is an isomorphism for each $n \geq 0$. 

While it may be possible to directly prove that the induced maps of derived functors are isomorphisms, it is desirable nonetheless to achieve a useful description of $\bb L_n(\fk m)$. For torsion-free ample groupoids we may describe $\bb L_n(\fk m)$ as the module-theoretic induced map in homology constructed in \cite{Miller24}. This is the culmination of work in \Cref{groupoid homology section} on the relationship between the module category of an ample groupoid and its equivariant Kasparov category.
\begin{thmx}[see \Cref{ample unit space compatibility} and \Cref{induced maps in homology agree}]\label{intro theorem on homology}
Let $G$ and $H$ be second countable Hausdorff ample groupoids with coefficient $\cs$-algebras $A \in \KK^G$ and $B \in \KK^H$, let $\Omega \colon G \to H$ be a second countable Hausdorff étale correspondence and let $f \in \KK^G(A,\Ind_\Omega B)$ be a morphism. Then the trivial families $\{G^0\}$ and $\{H^0\}$ are compatible under $\Omega$. For each $n \geq 0$ the map of derived functors $\bb L_n(\fk m)$ induced by the resulting ABC morphism $\fk m$ may be identified with the induced map in homology
\[ \H_n(\Omega;\K_*(f)) \colon \H_n(G;\K_*(A)) \to \H_n(H;\K_*(B)). \]
\end{thmx}

We arrive at \Cref{isomorphism theorem intro} through the combination of Theorems \ref{ABC thm}, \ref{intro theorem for obtaining topological K theory} and \ref{intro theorem on homology}. Outside of the torsion-free ample setting where $G$ and $H$ are merely étale and equipped with families $\cal E$ and $\cal F$ of open Morita embeddings, we may consider the associated derived functors $\bb L^{\cal E}_n \K_*(G \ltimes -)$ and $\bb L^{\cal F}_n \K_*(H \ltimes -)$ as exotic homology theories for $G$ and $H$. The establishment of more concrete frameworks for these homology theories could ultimately lead to new isomorphisms in K-theory.


The author would like to thank Christian Bönicke, Jamie Gabe, Xin Li and Valerio Proietti for discussions on the content of this paper, and Kevin Aguyar Brix, Jamie Gabe, Xin Li and Owen Tanner for detailed feedback on earlier drafts.

\section{The Baum--Connes conjecture for \'etale groupoids via localisation}\label{localisation section}

Throughout this paper, an étale groupoid is a topological groupoid with locally compact Hausdorff unit space whose range and source maps are local homeomorphisms. An ample groupoid is an étale groupoid with totally disconnected unit space.

Let $G$ be a second countable Hausdorff étale groupoid and let $A$ be a separable $G$-$\cs$-algebra. The Baum--Connes conjecture for $G$ with coefficients in $A$ asserts the bijectivity of the Baum--Connes assembly map 
\[ \mu_{G,A} \colon \Ktop_*(G;A) \to \K_*(G \ltimes_r A) \]
from the topological K-theory to the operator K-theory. Bönicke and Proietti use the Meyer--Nest framework of homological ideals in triangulated categories to describe the topological K-theory and the assembly map as a localisation. This is the picture of the Baum--Connes assembly map that we will use in this paper. This section serves primarily as a review of the relevant Meyer--Nest theory \cite{MeyNes06, MeyNes10, Meyer08} and its applications to étale groupoids and the Baum--Connes conjecture in \cite{ProYam22, BoePro24}.

We recall from e.g. \cite[Section 1.1]{BoePro24} the triangulated structure on the equivariant Kasparov category $\KK^G$ of separable $G$-$\cs$-algebras. The suspension functor $\Sigma \colon \KK^G \to \KK^G$ is given by $A \mapsto C_0((0,1),A)$ and satisfies $\Sigma^2 \cong \id$ by Bott periodicity. The \textit{mapping cone} of a $G$-equivariant $*$-homomorphism $f \colon A' \to B'$ is the $G$-$\cs$-algebra
\[ \Cone(f) = \{ (a,b) \in A' \oplus C_0((0,1],B') \mid f(a) = b(1) \}. \] 
The obvious $*$-homomorphisms $\Sigma B' \to \Cone(f)$ and $\Cone(f) \to A'$ together with Bott periodicity at $B'$ define the mapping cone triangle
\[ \Sigma B' \to \Cone(f) \to A' \to \Sigma^2 B'. \]
A triangle $A \to B \to C \to \Sigma A$ in $\KK^G$ is declared exact if it is isomorphic to a mapping cone triangle. Countable direct sums exist in $\KK^G$ and are given by the $\cs$-algebraic countable direct sums.

For any étale subgroupoid $H \subseteq G$, restriction defines a triangulated functor $\Res^H_G \colon \KK^G \to \KK^H$. Given an action $G \acts Z$ of $G$ on a locally compact Hausdorff space $Z$, any action of $G \ltimes Z$ on a $\cs$-algebra $A$ pushes forward to an action of $G$ on $A$. A $G$-$\cs$-algebra $A$ is \textit{proper} if there is a proper action $G \acts Z$ and an action $G \ltimes Z \acts A$ whose pushforward recovers $G \acts A$. 

In this paper a full subcategory of a triangulated category is \textit{localising} if it is a triangulated subcategory closed under countable direct sums. Localising subcategories are \textit{thick}, i.e. closed under direct summands. 

\begin{defn}
Let $G$ be a second countable Hausdorff étale groupoid. A $G$-$\cs$-algebra $A \in \KK^G$ is \textit{weakly contractible} if $\Res^H_G A = 0$ for each proper open subgroupoid $H \subseteq G$. We write $\fk N_G$ for the (localising) subcategory of $\KK^G$ of weakly contractible $G$-$\cs$-algebras and $\fk P_G$ for the localising subcategory generated by the proper $G$-$\cs$-algebras.
\end{defn}

Pairs of subcategories such as $(\fk P_G, \fk N_G)$ are key to (Bousfield) localisation in triangulated categories. 

\begin{defn}[Complementary pair of subcategories]
Let $\fk P$ and $\fk N$ be thick full triangulated subcategories of a triangulated category $\T$. We call the pair $(\fk P, \fk N)$ \textit{complementary} if
\begin{itemize}
\item $\T( P ,  N) = 0$ for each $P \in \fk P$, $N \in \fk N$, 
\item and for each $A \in \T$ there is an exact triangle $P \to A \to  N \to \Sigma P$ with $P \in \fk P$, $N \in \fk N$.
\end{itemize}
\end{defn}
We collect some basic results on complementary subcategories, see e.g. \cite[Proposition 2.9]{MeyNes06}:

\begin{prop}
Let $(\fk P, \fk N)$ be a complementary pair of subcategories of a triangulated category $\T$. 
\begin{itemize}
\item The categories determine each other: 
\begin{align*}
\fk N & = \{ N \in \T \mid \T(P,N) = 0 \; \text{for every} \; P \in \fk P \} \\
\fk P & = \{ P \in \T \mid \T(P,N) = 0 \; \text{for every} \; N \in \fk N \} 
\end{align*}
\item There is a functor $\Delta \colon A \mapsto \Delta_A = P_A \to A \to N_A \to \Sigma P_A$ taking an object $A$ of $\T$ to an exact triangle $\Delta_A$ in $\T$ with $P_A \in \fk P$ and $N_A \in \fk N$.
\item This determines $\Delta$ up to a unique natural isomorphism, and the components $A \mapsto P_A$ and $A \mapsto N_A$ are triangulated.
\item For $P \in \fk P$ and $N \in \fk N$, the maps $P_P \to P$ and $N \to N_N$ are isomorphisms.
\item Any triangulated functor $\Phi \colon  \T \to \T'$ vanishing on $\fk N$ factors uniquely up to natural isomorphism through the functor $A \mapsto P_A \colon  \T \to \fk P$.
\end{itemize}
\end{prop}

\begin{defn}\label{localisation definition}
With notation as above, the functor $A \mapsto P_A$ is known as the\footnote{Note that it is only defined up to unique natural isomorphism.} \textit{localisation functor} at $\fk N$, which we usually denote by $L \colon \T \to \T$. The morphism $P_A \to A$ for $A \in \T$ defines the \textit{localisation natural transformation} $\mu \colon L \Rightarrow \id$. For a functor $F \colon \T \to \fk C$, the functor $\bb L F \defeq F \circ L \colon \T \to \fk C$ is called the \textit{localisation} of $F$ at $\fk N$, and the natural transformation $\bb L F \Rightarrow F$ is the (localisation) \textit{assembly map}.
\end{defn}

Bönicke and Proietti prove that the Baum--Connes assembly map may be viewed as the localisation of the reduced operator K-theory functor at $\fk N_G$ \cite[Theorems B and C]{BoePro24}. 

\begin{thm}[Identification of the Baum--Connes assembly map]\label{baum connes assembly map identification}
Let $G$ be a second countable Hausdorff \'etale groupoid. Then $(\fk P_G, \fk N_G)$ is a complementary pair of subcategories of $\KK^G$ and for each $G$-$\cs$-algebra $A \in \KK^G$ the associated localisation assembly map 
\[\bb L \K_*(G \ltimes_r A) \to \K_*(G \ltimes_r A)\]
for the reduced K-theory functor $K_*(G \ltimes_r -)$ is naturally isomorphic to the Baum--Connes assembly map 
\[\mu_{G,A} \colon \Ktop_*(G;A) \to \K_*(G \ltimes_r A).\]
\end{thm}

\begin{rmk}\label{full reduced remark}
On the localising subcategory $\fk P_G$ generated by proper $G$-$\cs$-algebras, the reduced and full crossed products agree. The natural transformation from $G \ltimes -$ to $G \ltimes_r - $ applied to $LA \to A$ gives rise to a commutative diagram
\[\begin{tikzcd}
	{\bb L \K_*(G \ltimes A)} & {\K_*(G \ltimes A)} \\
	{\Ktop_*(G;A)} & {\K_*(G \ltimes_rA)}.
	\arrow[from=1-1, to=1-2]
	\arrow[from=1-2, to=2-2]
	\arrow["\cong"', from=1-1, to=2-1]
	\arrow["{\mu_{G,A}}"', from=2-1, to=2-2]
\end{tikzcd}\]
The Baum--Connes assembly map $\mu_{G,A}$ therefore factors through the canonical map $\K_*(G \ltimes A) \to \K_*(G \ltimes_r A)$, which is a surjection when $\mu_{G,A}$ is. It also means that the topological K-theory $\Ktop_*(G;A)$ enjoys the greater functoriality properties of the full crossed product.
\end{rmk}

To understand the pair $(\fk P_G, \fk N_G)$ we consider not just the objects but also morphisms which vanish under the restriction functors $\Res^H_G$ of proper open subgroupoids $H \subseteq G$. Such ideals are key to the Meyer--Nest framework for doing homological algebra in a triangulated category.

\begin{defn}[Homological functor and ideal]
Let $\T$ be a triangulated category with suspension functor $\Sigma$ and let $\fk C$ be an abelian category. A \textit{homological functor} $F \colon \T \to \fk C$ is an additive functor which maps exact triangles to long exact sequences. If $\fk C$ is stable with (additive) automorphism $\Sigma_{\fk C}$ a \textit{stable homological functor} is a homological functor $F \colon \T \to \fk C$ with a natural isomorphism $F \circ \Sigma \cong \Sigma_{\fk C} \circ F$. An ideal $\I \vartriangleleft \T$ of morphisms is a \textit{homological ideal} if it is the kernel of a stable homological functor.
\end{defn}
Typically our stable abelian category will just be the category $\Ab_*$ of $\Z$-graded abelian groups with the right shift automorphism. By the existence of a universal stable homological functor for any triangulated category, the kernel of any triangulated functor is homological (see \cite[Section 2.5]{MeyNes10}).

\begin{example}\label{homological ideal from a family}
Let $G$ be a second countable Hausdorff étale groupoid and let $\cal F$ be a family of proper open subgroupoids of $G$. Then 
\[ \I_{\cal F} \defeq \bigcap_{H \in \cal F} \ker \Res^H_G = \ker \prod_{H \in \cal F} \Res^H_G \]
is a homological ideal in $\KK^G$.
\end{example}

Many of the standard notions in homological algebra can be imitated relative to a homological ideal:

\begin{defn}
Let $\I$ be a homological ideal in a triangulated category $\T$, let $f \colon A \to B$ be a morphism in $\T$ embedded in an exact triangle
\[ \begin{tikzcd}
\Delta \arrow[r, "="{description}, draw=none] & A \arrow[r, "f"] & B \arrow[r, "g"] & C \arrow[r, "h"] & \Sigma A.
\end{tikzcd} \]
We say $f$ is \textit{$\I$-phantom} if $f \in \I$, \textit{$\I$-epic} if $g \in \I$, \textit{$\I$-monic} if $h \in \I$ and an \textit{$\I$-equivalence} if $g \in \I$ and $h \in \I$. An object is \textit{$\I$-contractible} if its identity morphism is in $\I$. An object $P \in \T$ is \textit{$\I$-projective} if $\I(P,A) = 0$ for every $A \in \T$, or equivalently for any morphism $f \colon P \to A$ and any $\I$-epimorphism $\pi \colon B \to A$ there is a lift $\tilde f \colon P \to B$ of $f$ through $\pi$. We say $\T$ has \textit{enough $\I$-projectives} or that $\I$ has enough projectives if for any object $A$ there is an $\I$-epimorphism $\pi \colon P \to A$ from an $\I$-projective object $P$. The exact triangle $\Delta$ is \textit{$\I$-exact} if $h \in \I$. A chain complex $C_\bullet = (C_n,d_n)_{n \in \Z}$ is \textit{$\I$-exact} if whenever we embed $d_n$ and $d_{n+1}$ in exact triangles 
\[ \begin{tikzcd}[column sep = 18]
C_n \arrow[r, "d_n"] & C_{n-1} \arrow[r, "f_n"] & X_n \arrow[r, "g_n"] & \Sigma C_n &  C_{n+1} \arrow[r, "d_{n+1}"] & C_{n} \arrow[r, "f_{n+1}"] & X_{n+1} \arrow[r, "g_{n+1}"] & \Sigma C_{n+1},
\end{tikzcd}  \]
we have $\Sigma f_{n+1} \circ g_n \in \I$. An \textit{$\I$-resolution} $P_\bullet \to A$ of an object $A$ is an $\I$-exact sequence $ \cdots \to P_n \to \cdots \to P_0 \to A \to 0$, which is an \textit{$\I$-projective resolution} if each $P_n$ is $\I$-projective.
\end{defn}

A homological ideal $\I \vartriangleleft \T$ can be used to construct a complementary pair of subcategories in $\T$. We denote by $\fk N_\I$ the (localising) full subcategory of $\I$-contractible objects $A \in \T$, and we write $\fk P_\I$ for the localising subcategory generated by the $\I$-projective objects in $\T$. For the homological ideal $\I_{\cal F}$ associated to a family $\cal F$ of open subgroupoids as in \Cref{homological ideal from a family} we may simply write $\fk N_{\cal F}$ and $\fk P_{\cal F}$.

\begin{thm}[Theorem 3.21 in \cite{Meyer08}]\label{complementary pair from an ideal}
Let $\T$ be a triangulated category with countable direct sums, and let $\I$ be a homological ideal of $\T$ compatible with countable direct sums and enough projectives. Then $(\fk P_\I , \fk N_\I)$ is a complementary pair of subcategories of $\T$.
\end{thm}

For clarity we may write $L_{\I} \colon \T \to \T $ for the associated localisation functor at $\fk N_\I$ and $\mu_{\I} \colon  L_{\I} \Rightarrow \id_{\T}\colon  \T \to \T$ for the localisation natural transformation. We write $\bb L^{\I} F$ for the localisation $F \circ L_\I$ of a functor $F$, but may write $\bb L F$ if the ideal is understood.

\begin{defn}[ABC tuple]\label{ABC tuple defn}
An \textit{ABC tuple} is a quadruple $(\T,\I,F,A)$ consisting of
\begin{itemize}
\item a triangulated category $\T$ with countable direct sums,
\item a homological ideal $\I \vartriangleleft \T$ closed under countable direct sums with enough projectives,
\item a stable homological functor $F \colon \T \to \Ab_*$ that commutes with countable direct sums,
\item and an object $A \in \T$.
\end{itemize}
\end{defn}
These are named after the ABC spectral sequence, covered at the end of this section, as an ABC tuple is precisely the input data for the construction of the ABC spectral sequence.

\begin{example}\label{subgroupoids ABC tuple}
Let $G$ be a second countable Hausdorff \'etale groupoid, let $A \in \KK^G$ be a $G$-$\cs$-algebra and let $F_G \colon \KK^G \to \Ab_*$ be the full K-theory functor $\K_*(G \ltimes -)$ or the reduced K-theory functor $\K_*(G \ltimes_r -)$. Let $\cal F$ be a countable family of open subgroupoids of $G$ with associated ideal $\I_{\cal F}$. Then 
\[ (\KK^G, \I_{\cal F}, F_G, A) \]
is an ABC tuple. As a special case, we may consider the trivial family $\cal F = \{ G^0 \}$, for which we write $\I_0$.
\end{example}

Let us explain why for countable $\cal F$ the homological ideal $\I_{\cal F}$ has enough projectives. For any open subgroupoid $H \subseteq G$ there is a triangulated induction functor $\Ind^G_H \colon \KK^H \to \KK^G$ and an adjunction $\Ind^G_H \dashv \Res^H_G$ \cite[Theorem 2.3]{BoePro24}. We therefore obtain an adjunction
\begin{align*}
\Ind_{\cal F} \colon \prod_{H \in \cal F} \KK^H & \to \KK^G  & \dashv & &\Res_{\cal F} \colon \KK^G & \to \prod_{H \in \cal F} \KK^H \\
(B_H)_{H \in \cal F} & \mapsto \bigoplus_{H \in \cal F} \Ind^G_H B_H && & A & \mapsto (\Res^H_G(A))_{H \in \cal F}.
\end{align*}
It is straightforward to verify that for any $H \in \cal F$ and $B \in \KK^H$ the induced algebra $\Ind^G_H B$ is $\I_{\cal F}$-projective, and that for any $A \in \KK^G$ the counit $\epsilon_A \colon \Ind_{\cal F} \Res_{\cal F} A \to A$ is $\I_{\cal F}$-epic. In general, an adjunction $E \dashv F$ of triangulated functors not only guarantees that $\I = \ker F$ has enough projectives but also provides a characterisation of the $\I$-projectives and an explicit $\I$-projective resolution of any object, see \cite[Section 3.6]{MeyNes10} and \cite[Section 2.1]{ProYam22}:
\begin{thm}\label{adjunctions are good}
Let $\T$ and $\T'$ be triangulated categories with countable direct sums and let $F \colon \T \to \T'$ be a triangulated functor with a left adjoint $E \colon \T' \to \T$ both compatible with countable direct sums. Then $\I = \ker F$ has enough projectives, and the $\I$-projective objects are precisely the direct summands of $EA$ for $A \in \T'$. Furthermore, for $B \in \T$ there is a natural $\I$-projective resolution
\begin{equation*}
\begin{tikzcd}
\cdots \arrow[r, "\delta_{n+1}"] & (EF)^{n+1}B \arrow[r, "\delta_n"] & \cdots \arrow[r, "\delta_2"] & (EF)^2B \arrow[r, "\delta_1"] & EFB \arrow[r, "\pi_0"] & B,
\end{tikzcd}
\end{equation*}
where $\pi_0 \defeq \epsilon_B \colon EFB \to B$ is given by the counit $\epsilon$ of the adjunction at $B$ and 
\begin{equation*}
\delta_n \defeq \sum_{i=0}^n (-EF)^i(\epsilon_{(EF)^{n-i} B}) \colon (EF)^{n+1}B \to (EF)^n B.
\end{equation*}
\end{thm}
From this description of the $\I$-projectives it follows that the localising subcategory $\langle E \rangle$ generated by the image of $E$ is equal to $\fk P_{\I}$ and complementary to $\fk N_\I$. For a countable family $\cal F$ of open subgroupoids as in \Cref{subgroupoids ABC tuple} the associated complementary pair is
\[(  \fk P_{\cal F} , \fk N_{\cal F} ) = (\langle \Ind_{\cal F} \rangle, \fk N_{\cal F}).\] 
Even when $\cal F$ is not countable, we write $\langle \Ind_{\cal F} \rangle $ for the localising subcategory generated by algebras induced by members of $\cal F$.
\begin{prop}
Let $G$ be a second countable Hausdorff étale groupoid and let $\cal F$ be a (not necessarily countable) family of proper open subgroupoids of $G$. Then we have containments
\begin{align*}
 \langle \Ind_{\cal F} \rangle & \subseteq \fk P_G, & \fk N_G & \subseteq \fk N_{\cal F}. 
\end{align*}
\begin{proof}
For the first containment, for each $P \in \cal F$ and $B \in \KK^P$ the induced algebra $\Ind^G_P B$ is a proper $G$-$\cs$-algebra through the proper $G$-space $G/P$. The containment $\fk N_G \subseteq \fk N_{\cal F}$ is immediate. 
\end{proof}
\end{prop}

In order to prove \Cref{baum connes assembly map identification}, Bönicke and Proietti construct a countable family $\cal F$ of proper open subgroupoids of $G$ for which the reverse inclusions also hold \cite[Theorem 3.10 and Lemma 3.16]{BoePro24}, and so
\[ (\fk P_G, \fk N_G) = (\fk P_{\cal F} , \fk N_{\cal F}). \]
Diving into the details of the construction, the family $\cal F$ may moreover be chosen to be finer than any pre-specified open cover $\cal U$ of $G^0$ in the sense that for any $H \in \cal F$ there is $U \in \cal U$ with $H^0 \subseteq U$. As a consequence, we obtain the following:
\begin{cor}\label{open cover vanishing}
Let $G$ be a proper second countable Hausdorff \'etale groupoid and let $\cal U$ be an open cover of $G^0$. Then $A \in \KK^G$ vanishes if and only if $\Res^{G \restriction_U}_G A$ vanishes for each $U \in \cal U$.
\begin{proof}
Suppose that $\Res^{G \restriction_U}_G A \cong 0$ for each $U \in \cal U$. Choose a countable family $\cal I$ of proper open subgroupoids with $\fk N_{\cal I} = \fk N_G$ such that for any $H \in \cal I$ there is $U \in \cal U$ with $H^0 \subseteq U$. It follows that $\Res^H_G A \cong 0$ for each $H \in \cal I$, and so $A \in \fk N_G$. As $G$ itself is a proper open subgroupoid of $G$, we have $A \cong \Res^G_G A \cong 0$.
\end{proof}
\end{cor}
Framing the Baum--Connes assembly map in terms of localisation gives us access to the ABC spectral sequence, which converges to the localisation. This was utilised by Proietti and Yamashita to obtain a spectral sequence \cite[Theorem 4.3]{ProYam22}
\[ E^2_{p,q} = \H_p(G;K_q(A)) \Rightarrow \K_{p+q}(G \ltimes_r A) \]
relating groupoid homology to operator K-theory for a second countable Hausdorff ample groupoid $G$ and a separable $G$-$\cs$-algebra $A$, where $G$ has torsion-free isotropy groups and satisfies a strengthening of the Baum--Connes conjecture. On the second sheet of the general ABC spectral sequence are the derived functors with respect to a homological ideal:
\begin{defn}[Derived functors]
Let $\I$ be a homological ideal in a triangulated category $\T$ with enough projectives and let $F \colon \T \to \Ab$ be a functor. For each $n \geq 0$, the $n$th \textit{left derived functor} $\bb L_n^{\I} F \colon \T \to \Ab$ of $F$ is given as follows. For $A \in \T$ pick any $\I$-projective resolution
\[ \cdots \to P_n \to \cdots \to P_0 \to A  \]
of $A$. The group $\bb L_n^{\I} F (A)$ is the $n$th homology group $\H_n(F(P_\bullet))$ of the chain complex
\[ \cdots \to F(P_n) \to \cdots \to F(P_0) \to 0. \]
For $n \leq -1$ we set $\bb L^{\I}_n F = 0$. 
\end{defn}
When the homological ideal is clear from context we may simply write $\bb L_n$, and for the homological ideal $\I_{\cal F}$ associated to a family $\cal F$ of subgroupoids we write $\bb L^{\cal F}_n$ instead of $\bb L^{\I_{\cal F}}_n$. Independence from the choice of $\I$-projective resolution and the functoriality of this construction follow from a version of the fundamental lemma of homological algebra for $\I$-projective resolutions \cite[Proposition 3.26]{MeyNes10}:
\begin{lem}[The fundamental lemma of homological algebra]\label{resolutionmap}
Let $\T$ be a triangulated category and $\I$ be a homological ideal. If $A$ and $B$ are objects of $\T$, $P_\bullet \rightarrow A$ is an $\I$-projective resolution of $A$ and $Q_\bullet \rightarrow B$ is an $\I$-resolution of $B$, then any morphism $f\colon A \rightarrow B$ lifts to a chain map $\tilde{f}\colon P_\bullet \rightarrow Q_\bullet$ uniquely up to chain homotopy. 
\end{lem}

\begin{example}[Groupoid homology]\label{groupoid homology example}
Let $G$ be a second countable Hausdorff ample groupoid and consider the ideal $\I_0 \vartriangleleft \KK^G$ for the trivial family $\{G^0\}$. Then the $p$th derived functor of the K-theory functor $\K_q(G \ltimes_r -)$ at a $G$-$\cs$-algebra $A \in \KK^G$ is naturally isomorphic to the groupoid homology $\H_p(G;\K_q(A))$ of $G$ with coefficients in the $G$-module $\K_q(A)$ \cite[Proposition 4.2]{ProYam22}:
\[ \bb L^{\I_0}_p \K_q(G \ltimes_r A) \cong \H_p(G;\K_q(A)) \]
This also holds for the full crossed product as it agrees with the reduced crossed product on $\I_0$-projective objects.
\end{example}

The ABC spectral sequence links together the derived functors $\bb L^\I_n F(A)$ and the localisation $\bb L^\I F(A)$ \cite[Theorem 5.1]{Meyer08}:

\begin{thm}[The ABC spectral sequence]\label{ABC}
Let $(\T,\I,F,A)$ be an ABC tuple. Then there is a strongly convergent spectral sequence 
\[ E^2_{p,q} = \bb{L}^\I_p F_q (A) \Rightarrow \bb{L}^\I F_{p+q}(A). \]
\end{thm}

\begin{example}\label{ABC spectral sequence for top K theory}
Let $G$ be a second countable Hausdorff \'etale groupoid, let $\cal F$ be a countable family of proper open subgroupoids of $G$ and let $A \in \KK^G$ be a $G$-$\cs$-algebra. If $\fk N_{\cal F} = \fk N_G$ (and therefore also $\langle \Ind_{\cal F} \rangle = \fk P_G$) then the ABC spectral sequence for $(\KK^G,\I_{\cal F}, \K_*(G \ltimes -),A)$ is given by
\[ \bb L^{\cal F}_p \K_q(G \ltimes A) \Rightarrow \Ktop_{p+q}(G;A). \] 
\end{example}

\section{Homological ideals associated to families of openly embedded groupoids}\label{family study}

We have seen that we may study the topological K-theory $\Ktop_*(G;A)$ of an étale groupoid $G$ with coefficients in a $G$-$\cs$-algebra $A$ through a family $\cal F$ of proper open subgroupoids of $G$ satisfying
\[ (\langle \Ind_{\cal F} \rangle, \fk N_{\cal F}) =  (\fk P_G, \fk N_G).  \]
For each family with this property we obtain an instance of the ABC spectral sequence that converges to the topological K-theory (\Cref{ABC spectral sequence for top K theory}). This raises questions of how to construct families with this property, how to check a family has this property and how families with this property relate to each other. We introduce the following sufficient condition:

\begin{defn}[Condition (P)]
Let $G$ be an \'etale groupoid with a family $\cal F$ of proper open subgroupoids. We say that $\cal F$ satisfies \textit{condition (P)} if for any proper open subgroupoid $P \subseteq G$ and any $x \in P^0$ there is an open neighbourhood $U \subseteq P^0$ of $x$ and a member $H \in \cal F$ such that $P \restriction_U \subseteq H$.
\end{defn}

\begin{prop}\label{condition P gives topological K theory}
Let $G$ be a second countable Hausdorff \'etale groupoid and let $\cal F$ be a countable family of proper open subgroupoids satisfying condition (P). Then 
\[ (\langle \Ind_{\cal F} \rangle, \fk N_{\cal F}) =  (\fk P_G, \fk N_G).  \]
\begin{proof}
Let $A \in \fk N_{\cal F}$ and let $P \subseteq G$ be a proper open subgroupoid. Applying \Cref{open cover vanishing} we have $\Res^P_G A = 0$ by construction and so $A \in \fk N_G$. As complementary subcategories determine each other, we also get $\langle \Ind_{\cal F} \rangle = \fk P_G$.
\end{proof}
\end{prop}

To check condition (P) in practice we need only check that for each finite subgroup $\Gamma$ of the isotropy of $G$ there is a member of $\cal F$ containing $\Gamma$. This is due to the following simple observation.

\begin{lem}\label{locally containing proper groupoids observation}
Let $P$ be a proper \'etale groupoid and let $O \subseteq P$ be an open subset containing the isotropy $P^x_x$ at some point $x \in P^0$. Then $O$ contains $P \restriction_U$ for some open neighbourhood $U \subseteq P^0$ of $x$.
\begin{proof}
The image of $P \setminus O$ under the closed map $r \times s \colon P \to P^0 \times P^0$ is disjoint from $(x,x)$. It therefore misses $U \times U$ for an open neighbourhood $U$ of $x$, and consequently $P \restriction_U \subseteq O$.
\end{proof}
\end{lem}

Conversely, each finite subgroup of isotropy of $G$ is contained in some proper open subgroupoid.

\begin{lem}\label{proper subgroupoid containing isotropy}
Let $S \acts X$ be an action of an inverse semigroup $S$ by partial homeomorphisms on a locally compact Hausdorff space $X$. For any $x \in X$ and finite subgroup $\Gamma \leq (S \ltimes X)^x_x$ there is an idempotent $e \in E(S)$ with $x \in \dom e$ and a finite subgroup $F \leq S_e = \{ s \in S \suchthat s^*s = s s^* = e \}$ such that $\Gamma = \{ [s,x] \suchthat s \in F \}$. Consequently, $\Gamma$ is contained in the proper open subgroupoid $\{[s,y] \suchthat s \in F, \; y \in \dom e \}$.
\begin{proof}
For each $g \in \Gamma$, let $s_g \in S$ represent $g$ in the sense that $g = [s_g,x]$, with an idempotent representing the identity element of $\Gamma$. For each pair $g,h \in \Gamma$ let $e_{g,h} \in E(S)$ be an idempotent with $x \in \dom e_{g,h}$, such that $s_g s_h e_{g,h} = s_{gh} e_{g,h}$, $e_{g,h} \leq \dom s_h$, $e_{g,h} \leq \dom s_{gh}$ and $s_h \bcdot e_{g,h} \leq \dom s_g$. Let $e_0 = \prod_{g,h \in \Gamma} e_{g,h}$ and let $e = \prod_{g \in \Gamma} s_g \bcdot e_0$. For each $g \in \Gamma$, the element $t_g = s_g e$ belongs to $S_e$ because $s_g \bcdot e = e$. The finite set $F = \{t_g \suchthat g \in \Gamma \}$ is a subgroup of $S_e$ and by construction satisfies $\Gamma = \{[t,x] \suchthat t \in F \}$.
\end{proof}
\end{lem}

The following characterisation of condition (P) follows straightforwardly from the previous two lemmas, noting that every \'etale groupoid may be written as the transformation groupoid of an inverse semigroup action.
\begin{prop}\label{checkable P}
Let $G$ be an \'etale groupoid and let $\cal F$ be a family of proper open subgroupoids of $G$. Then $\cal F$ satisfies condition (P) if and only if for each $x \in G^0$ and finite subgroup $\Gamma \leq G^x_x$, there is a member $K \in \cal F$ with $\Gamma \subseteq K$.
\end{prop}

The construction of proper open subgroupoids containing finite isotropy subgroups in \Cref{proper subgroupoid containing isotropy} gives us a recipe for constructing families of subgroupoids satisfying condition (P).

\begin{example}\label{inverse semigroup P example}
Let $G = S \ltimes X$ be the transformation groupoid of an action $S \acts X$ of an inverse semigroup $S$ by partial homeomorphisms on a locally compact Hausdorff space $X$. For each non-zero idempotent $e \in E(S)^\times$ and finite subgroup $F \leq S_e$ consider the proper open subgroupoid $P(e,F) = \{[s,y] \suchthat s \in F, \; y \in \dom e \}$ of $G$. Then the family 
\[ \cal F_{\text{fin}} = \{ P(e,F) \suchthat e \in E(S)^\times, \; F \leq S_e \; \text{finite subgroup}  \}  \]
of proper open subgroupoids of $G$ satisfies condition (P) and is countable when $S$ is countable. The subgroupoid $P(e,F)$ may be identified with the transformation groupoid $F_\downarrow \ltimes X$, where $F_\downarrow = \{ f d \mid d \in E(S), f \in F \}$ is the inverse semigroup of restrictions of elements of $F$ in $S$.
\end{example}

\begin{rmk}
For a second countable étale groupoid $G$ the existence of a countable family $\cal E$ of subgroupoids satisfying (P) follows from the above example. It follows that any other family $\cal F$ satisfying (P) has a countable subfamily satisfying (P), since we need only check the condition against the countable family $\cal E$, each member of which needs only countably many members of $\cal F$ by second countability.
\end{rmk}

Suppose $G$ and $H$ are Morita equivalent étale groupoids and let $\cal F$ be a countable family of proper open subgroupoids of $G$. The Morita equivalence induces an equivalence $\KK^G \cong \KK^H$ of Kasparov categories, so we may search for a counterpart in $\KK^H$ to the homological ideal $\I_{\cal F} \vartriangleleft \KK^G$. However, there is no clear way to transfer subgroupoids of $G$ to subgroupoids of $H$. For example, there is no subgroupoid of the trivial groupoid $R_1$ which under the Morita equivalence with the full equivalence relation $R_2$ on two points corresponds to the unit space $R_2^0 \subseteq R_2$. Instead we consider a looser notion of subgroupoid which is morally an open subgroupoid up to Morita equivalence, or an open embedding in a category of étale groupoids with étale correspondences as morphisms. We refer the reader to \cite{AKM22} and \cite[Section 3.1]{Miller23a} for an introduction to and examples of \'etale correspondences.

\begin{defn}[Étale correspondence]
Let $G$ and $H$ be \'etale groupoids. An \textit{\'etale correspondence} $\Omega \colon G \to H$ is a $G$-$H$-bispace $\Omega$ with anchor maps $\rho\colon \Omega \to G^0$ and $\sigma \colon  \Omega \to H^0$ called the \textit{range} and \textit{source} such that the right action $\Omega \rightacts H$ is free, proper and \'etale. An étale correspondence $\Omega \colon G \to H$ is \textit{proper} if the induced map $\orho \colon \Omega/H \to G^0$ is proper and \textit{tight} if $\orho$ is an homeomorphism. The \textit{composition} $\Lambda \circ \Omega$ of \'etale correspondences $\Omega \colon G \to H$ and $\Lambda \colon H \to K$ is given by the fibre product $\Omega \times_H \Lambda$ over $H$. A \textit{slice} $U \subseteq \Omega$ is an open set on which the source $\sigma \colon U \to H^0$ and the orbit map $q \colon U \to \Omega/H$ are injective.
\end{defn}
Composition of étale correspondences is associative up to isomorphism, with an identity up to isomorphism at $G$ given by the $G$-$G$-bispace $G$. Invertible correspondences are precisely Morita equivalences.

\begin{defn}[Open Morita embedding]
Let $G$ and $P$ be \'etale groupoids. An \textit{open Morita embedding} $P \hookrightarrow_M G$ of $P$ into $G$ is a tight \'etale correspondence $\Omega \colon P \to G$ with a free left action $P \acts \Omega$. We say that $P$ is \textit{openly Morita embedded} into $G$.
\end{defn}

The correspondence $G^{P^0} \colon P \to G$ associated to an open inclusion $P \subseteq G$ is an open Morita embedding, as is a Morita equivalence. The composition of open Morita embeddings is again an open Morita embedding. These compositions describe all open Morita embeddings thanks to the following construction which appears in \cite[Proposition 7.16]{MeyZhu15}. For an \'etale groupoid $H$ and a free, proper \'etale right $H$-space $\Omega$, we define an \'etale groupoid structure on the space $\Omega \times_H \Omega^*$. Here $\Omega^*$ is a formal \textit{opposite} topological space to $\Omega$ defined by a homeomorphism
\[ \omega \mapsto \omega^* \colon \Omega \to \Omega^*, \]
equipped with a left action $H \acts \Omega^*$ determined by $h \bcdot \omega^* \defeq (\omega \bcdot h^{-1})^*$ for $\omega \in \Omega$ and $h \in H_{\sigma(\omega)}$. Define range, source, composition and inversion maps on $\Omega \times_H \Omega^*$ by
\begin{align*}
r([\omega_1, \omega_2^*]_H) & = [\omega_1, \omega_1^*]_H \\
s([\omega_1, \omega_2^*]_H) & = [\omega_2, \omega_2^*]_H \\
[\omega_1, \omega_2^*]_H [\omega_2, \omega_3^*]_H & = [\omega_1, \omega_3^*]_H \\
[\omega_1, \omega_2^*]_H^{-1} & = [\omega_2, \omega_1^*]_H
\end{align*}
for $\omega_1, \omega_2, \omega_3 \in \Omega$ with $\sigma(\omega_1) = \sigma(\omega_2) = \sigma(\omega_3)$. This defines an \'etale groupoid whose unit space may be identified with $\Omega/H$. Furthermore, $\Omega$ serves as a Morita equivalence $\Omega \times_H \Omega^* \sim_M H \restriction_{\sigma(\Omega)}$ through an action $\Omega \times_H \Omega^* \acts \Omega$ determined by 
\[[\omega_1, \omega_2^*]_H \bcdot \omega_2 = \omega_1\]
for $\omega_1, \omega_2 \in \Omega$ with $\sigma(\omega_1) = \sigma(\omega_2)$. We note that if $\Omega$ is Hausdorff or second countable then so is $\Omega \times_H \Omega^*$.

\begin{prop}\label{decomposition of correspondences}
Let $\Omega \colon G \to H$ be an \'etale correspondence. Then there is an \'etale groupoid $K$ and a decomposition of $\Omega$ into 
\begin{itemize}
\item an actor $G \acts \Omega \times_H \Omega^*$,
\item the Morita equivalence $\Omega \colon \Omega \times_H \Omega^* \to H \restriction_{\sigma(\Omega)}$, and
\item the open inclusion $H\restriction_{\sigma(\Omega)} \subseteq H$.
\end{itemize}
If $\Omega$ is an open Morita embedding then the actor $G \acts \Omega \times_H \Omega^*$ arises from an open embedding which is a homeomorphism on unit spaces. 
\begin{proof}

Define an actor $G \acts \Omega \times_H \Omega^*$ by $g \bcdot [\omega_1, \omega_2^*]_H \defeq [g \bcdot \omega_1, \omega_2^*]_H$ for $[\omega_1, \omega_2^*]_H \in \Omega \times_H \Omega^*$ and $g \in G_{\rho(\omega_1)}$. The composition of the actor $G \acts \Omega \times_H \Omega^*$, the Morita equivalence $\Omega \times_H \Omega^* \sim_M H\restriction_{\sigma(\Omega)}$ and the open inclusion $H\restriction_{\sigma(\Omega)} \subseteq H$ is isomorphic to $\Omega \colon G \to H$ via the isomorphism
\begin{align*} 
(\Omega \times_H \Omega^*) \times_{\Omega \times_H \Omega^*} \Omega \times_{H\restriction_{\sigma(\Omega)}} H\restriction_{\sigma(\Omega)} & \to \Omega  \\
[[\omega_1,\omega_2^*],\omega_2,h] & \mapsto \omega_1 \bcdot h
\end{align*}
of $G$-$H$-bispaces, where $\omega_1, \omega_2 \in \Omega$ and $h \in H$ satisfy $\sigma(\omega_1) = \sigma(\omega_2) = r(h)$.

If $\Omega \colon G \to H$ is an open Morita embedding, then so is the actor $G \acts \Omega \times_H \Omega^*$. By tightness the left anchor map induces a homeomorphism of unit spaces $\varphi \colon G^0 \to \Omega/H$. The continuous homomorphism $g \mapsto g \bcdot \varphi(s(g)) \colon G \to \Omega \times_H \Omega^*$ is \'etale as it restricts to $\varphi$ on the unit spaces, and it is injective by freeness of $G \acts \Omega \times_H \Omega^*$.
\end{proof}
\end{prop}

This also tells us that an open Morita embedding arises from an open embedding if the Morita equivalence $\Omega \times_H \Omega^* \sim_M H \restriction_{\sigma(\Omega)}$ is an isomorphism. This happens precisely if $\Omega = U \bcdot H$ for a slice $U \subseteq \Omega$, in which case we may identify $U$ with $\sigma(\Omega)$ and $\Omega$ with $H \restriction_{\sigma(\Omega)}$. However, even the open Morita embeddings which are not open embeddings resemble embeddings locally. In general, tight \'etale correspondences locally resemble groupoid homomorphisms which are open embeddings on the level of unit spaces:

\begin{prop}\label{recognising embeddings}
Let $\Omega \colon G \to H$ be a tight \'etale correspondence. Then for each slice $U \subseteq \Omega$ there is an \'etale homomorphism $\varphi_U \colon G \restriction_{\rho(U)} \to H$ such that $\varphi_U(\rho(u)) = \sigma(u)$ for each $u \in U$ and $g \bcdot u = u \bcdot \varphi_U(g)$ for each pair $(g,u) \in G \times_{G^0} U$. This \'etale homomorphism induces the composition
\[ G \restriction_{\rho(U)} \xhookrightarrow{\iota} G \xrightarrow{\Omega} H  \]
on the level of correspondences. When $\Omega \colon G \to H$ is further an open Morita embedding, then $\varphi_U$ is an open embedding. In particular, if $U \bcdot H = \Omega$ then $\Omega$ is induced by an open embedding.
\begin{proof}
Let $\rho_U \colon U \to G^0$ be the restriction to $U$ of the anchor map for the action $G \acts \Omega$ and consider the following continuous maps.
\begin{align*}
G \restriction_{\rho(U)} & \to \Omega^* \times_{\Omega/H} \Omega & \Omega^* \times_{\Omega/H} \Omega & \to H \\
g & \mapsto (\rho_U^{-1}(r(g))^*, g \bcdot \rho_U^{-1}(s(g))) & (\omega^*, \omega \bcdot h) & \mapsto h 
\end{align*}
The first is well-defined by tightness of $\Omega$ and the second is well-defined and continuous because $\Omega \rightacts H$ is free, proper and \'etale. We set $\varphi_U \colon G \restriction_{\rho(U)} \to H$ to be the composition of these maps. By construction this is a homomorphism and satisfies $\varphi_U(\rho(u)) = \sigma(u)$ for each $u \in U$ and $g \bcdot u = u \bcdot \varphi_U(g)$ for each pair $(g,u) \in G \times_{G^0} U$. The correspondence associated to $\varphi_U$ has bispace $\rho(U) \times_{H^0} H$, while the composition $\Omega \circ \iota$ has bispace $G^{\rho(U)} \times_G \Omega$. The homeomorphism
\begin{align*}
\rho(U) \times_{H^0} H & \to G^{\rho(U)} \times_G \Omega \\
(x,h) & \mapsto [x, \rho_U^{-1}(x) \bcdot h]_G
\end{align*}
is an isomorphism of correspondences. Note that these spaces are both homeomorphic to $U \bcdot H \subseteq \Omega$. 

If we now assume that $\Omega \colon G \to H$ is an open Morita embedding then so is the correspondence associated to $\varphi_U$, which implies that $\varphi_U$ is an \'etale injection and thus an open embedding. When we further have $U \bcdot H = \Omega$ then by tightness $\rho(U) = G^0$ and so $\varphi_U$ induces $\Omega$.
\end{proof}
\end{prop}

To define a restriction functor for an openly Morita embedded groupoid we utilise the induction functor $\Ind_\Omega \colon \KK^H \to \KK^G$ associated to a second countable Hausdorff étale correspondence $\Omega \colon G \to H$ constructed in \cite[Section 6]{Miller23a}. This is a triangulated functor because it sends mapping cones of $H$-equivariant $*$-homomorphisms to mapping cones of $G$-equivariant $*$-homomorphisms. The induction functor of the étale correspondence $H^{G^0} \colon G \to H$ associated to an open inclusion $G \subseteq H$ is naturally isomorphic to the restriction functor $\Res^G_H$. This motivates the following definition:
\begin{defn}[Restriction for open Morita embeddings]
Let $G$ and $P$ be second countable Hausdorff étale groupoids and let $P \hookrightarrow_M G$ be an open Morita embedding. The \textit{restriction functor} $\Res^P_G \colon \KK^G \to \KK^P$ is the induction functor of the étale correspondence $P \hookrightarrow_M G$.
\end{defn}
The induction functors associated to Morita equivalences are themselves equivalences. As an open Morita embedding $P \hookrightarrow_M G$ is the composition of open inclusions and Morita equivalences, the restriction functor $\Res^P_G \colon \KK^G \to \KK^P$ has a left adjoint which we term the induction functor $\Ind^G_P \colon \KK^P \to \KK^G$. Thus by \Cref{adjunctions are good}, for a countable family $\cal F$ of second countable Hausdorff groupoids openly Morita embedded into $G$, the localising subcategory $\langle \Ind_{\cal F} \rangle$ generated by these induced algebras is complementary to the full subcategory $\fk N_{\cal F}$ of $G$-$\cs$-algebras $B$ for which $\Res^H_G B \cong 0$ for each $H \in \cal F$. This enables us to build ABC tuples from families of open Morita embeddings, extending \Cref{subgroupoids ABC tuple}:

\begin{example}\label{OMEs ABC tuple}
Let $G$ be a second countable Hausdorff \'etale groupoid, let $A \in \KK^G$ be a $G$-$\cs$-algebra and let $F_G \colon \KK^G \to \Ab_*$ be the full K-theory functor $\K_*(G \ltimes -)$ or the reduced K-theory functor $\K_*(G \ltimes_r -)$. Let $\cal F$ be a countable family of second countable Hausdorff groupoids openly Morita embedded into $G$ and set
\[ \I_{\cal F} = \bigcap_{H \in \cal F} \ker \Res^H_G \vartriangleleft \KK^G. \]
Then $(\KK^G, \I_{\cal F}, F_G, A)$ is an ABC tuple.
\end{example}

\begin{prop}
Let $G$ be a second countable Hausdorff étale groupoid and let $\cal F$ be a family of proper second countable Hausdorff groupoids openly Morita embedded into $G$. Then 
\begin{align*}
 \langle \Ind_{\cal F} \rangle & \subseteq \fk P_G, & \fk N_G & \subseteq \fk N_{\cal F}. 
\end{align*}
\begin{proof}
It suffices to consider singleton families. In this case $(\langle \Ind_{\cal F} \rangle, \fk N_{\cal F})$ is complementary, so we need only demonstrate one inclusion. Let $\Omega \colon P \hookrightarrow_M G$ be an open Morita embedding with $P$ proper and let $B \in \KK^P$. Recall that the induction functor $\Ind^G_P$ is a composition of induction functors for the decomposition of $P \hookrightarrow_M G$ into open inclusions and Morita equivalences (\Cref{decomposition of correspondences}):
\[ P \hookrightarrow \Omega \times_H \Omega^* \sim_M G \restriction_{\sigma(\Omega)} \hookrightarrow G \]
For both of these open inclusions the included groupoid acts properly on the including groupoid. For open inclusions $H \subseteq K$ with this property, the induced algebra $\Ind^K_H A$ of an $H$-$\cs$-algebra $A$ is $\KK^K$-equivalent to the generalised fixed-point algebra $s^* A^H$, where $s \colon K_{H^0} \to H^0$ is the source map \cite[Remark 2.1]{BoePro24}. The fixed point algebra $s^* A^H$ is none other than the induced algebra $\Ind_\Lambda A$ for the étale correspondence $\Lambda = K_{H^0} \colon K \to H$. It follows that $\Ind^G_P B \cong \Ind_{\Omega^*} B$ in $\KK^G$ for any $B \in \KK^P$, where $\Omega^* \colon G \to P$ is the étale correspondence given by the opposite bispace. This correspondence factors through the action correspondence associated to the proper $G$-space $\Omega^*/P$, which induces the forgetul functor $\KK^{G \ltimes \Omega^*/P} \to \KK^G$. Therefore $\Ind^G_P B$ is $\KK^G$-equivalent to a proper $G$-$\cs$-algebra, and so $\Ind^G_P B \in \fk P_G$.
\end{proof}
\end{prop}

\begin{defn}
Let $G$ be a second countable Hausdorff \'etale groupoid with a family $\cal F$ of proper openly Morita embedded groupoids. We say that $\cal F$ satisfies \textit{condition (P)} if for any proper open subgroupoid $P \subseteq G$ and any $x \in P^0$ there is an open neighbourhood $U \subseteq P^0$ of $x$ and a member $H \in \cal F$ such that $P\restriction_U \hookrightarrow G$ factors through $H \hookrightarrow_M G$.
\end{defn}

It is straightforward to see that this agrees with the previous definition in the case that $\cal F$ consists entirely of genuine subgroupoids. This condition can also be checked at each finite isotropy subgroup:

\begin{prop}\label{checkable general P}
Let $G$ be an \'etale groupoid with a family $\cal F$ of proper openly Morita embedded groupoids. Then $\cal F$ satisfies condition (P) if and only if for each $x \in G^0$ and finite subgroup $\Gamma \leq G^x_x$, there is a member $H \in \cal F$ with embedding $\Omega \colon H \hookrightarrow_M G$ and an element $\omega \in \Omega$ with $\sigma(\omega) = x$ and $\omega \bcdot \Gamma \subseteq H \bcdot \omega$.
\begin{proof}
Suppose that $\cal F$ satisfies condition (P) and let $x \in G^0$ and $\Gamma \subseteq G^x_x$ be as in the statement. Then by \Cref{proper subgroupoid containing isotropy} there is a proper open subgroupoid $P \subseteq G$ containing $\Gamma$. By hypothesis there is an open neighbourhood $U \subseteq P^0$ of $x$ and a member $H \in \cal F$ such that $P\restriction_U \to G$ factors through the embedding $\Omega \colon H \hookrightarrow_M G$ via an \'etale correspondence $\Lambda \colon P\restriction_U \to H$. We obtain an isomorphism $\Phi \colon G^U \cong \Lambda \times_H \Omega$ of $P\restriction_U$-$G$-bispaces through which $\Phi(x) = [\lambda,\omega]_H$ for some $\lambda \in \Lambda$ and $\omega \in \Omega$. For each $g \in \Gamma$ we may calculate that
\[ [\lambda, \omega \bcdot g]_H = \Phi(x \bcdot g) = \Phi(g \bcdot x) = [g \bcdot \lambda, \omega]_H \]
and consequently there is $h \in H$ such that $h \bcdot \omega = \omega \bcdot g$.

For the converse, let $P \subseteq G$ be a proper open subgroupoid and let $x \in P^0$. Then $\Gamma = P^x_x$ is finite so by hypothesis we obtain a member $H \in \cal F$ with embedding $\Omega \colon H \hookrightarrow_M G$ and $\omega \in \Omega$ with $\sigma(\omega) = x$ and $\omega \bcdot \Gamma \subseteq H \bcdot \omega$. Let $V \subseteq \Omega$ be a slice containing $\omega$. By \Cref{recognising embeddings} the composition of $\Omega \colon H \hookrightarrow_M G$ with the inclusion $H \restriction_{\rho(V)} \hookrightarrow H$ is induced by an open embedding $\varphi_V \colon H \restriction_{\rho(V)} \hookrightarrow G$ satisfying $\varphi_V(\rho(v)) = \sigma(v)$ for each $v \in V$ and $h \bcdot v = v \bcdot \varphi_V(h)$ for each pair $(h,v) \in H \times_{H^0} V$. For each $g \in \Gamma$ pick $h \in H$ such that $\omega \bcdot g = h \bcdot \omega$, and thus $h \in H^{\rho(\omega)}_{\rho(\omega)}$. By construction $\omega \bcdot \varphi_V(h) = \omega \bcdot g$ and so $g = \varphi_V(h)$ is in the image of $\varphi_V$. This image $\im \varphi_V \subseteq G$ is open and contains $\Gamma$ and thus by \Cref{locally containing proper groupoids observation} contains the restriction of $P$ to some open neighbourhood $U$ of $x$. By inverting we obtain an open embedding $\varphi_V^{-1} \colon P \restriction_U \hookrightarrow H$ which composes with the open Morita embedding $H \hookrightarrow_M G$ to recover the inclusion $P \restriction_U \hookrightarrow G$.
\end{proof}
\end{prop}

\begin{prop}\label{properly weak equivalence result}
Let $G$ be a second countable Hausdorff \'etale groupoid with a family $\cal F$ of proper second countable Hausdorff openly Morita embedded groupoids satisfying condition (P). Then 
\[ (\langle \Ind_{\cal F} \rangle, \fk N_{\cal F}) =  (\fk P_G, \fk N_G).  \]
\begin{proof}
We may without loss of generality assume that $\cal F$ is countable as we may pass to a countable subfamily $\cal E$ which also satisfies (P). The pair $(\langle \Ind_{\cal F} \rangle, \fk N_{\cal F})$ is then complementary, so it suffices to prove that $\fk N_{\cal F} = \fk N_G$. Let $A \in \fk N_{\cal F}$, let $P \subseteq G$ be a proper open subgroupoid and let $x \in P^0$. Pick an open neighbourhood $U \subseteq P^0$ of $x$ and a member $Q \in \cal F$ such that $P \restriction_U \hookrightarrow G$ factors through $Q \hookrightarrow_M G$. It follows that $\Res^{P \restriction_U}_G A \cong 0$. Such $U$ cover $P^0$, so by \Cref{open cover vanishing} we may conclude that $\Res^P_G A \cong 0$. 
\end{proof}
\end{prop}

\begin{example}
Let $G$ be a second countable Hausdorff \'etale groupoid with a countable family $\cal F$ of proper openly Morita embedded groupoids and let $A \in \KK^G$ be a $G$-$\cs$-algebra. If $\cal F$ satisfies condition (P) then the ABC spectral sequence for $(\KK^G,\I_{\cal F}, \K_*(G \ltimes -),A)$ is given by
\[ \bb L^{\cal F}_p \K_q(G \ltimes A) \Rightarrow \Ktop_{p+q}(G;A). \] 
If $G$ has torsion-free stabilisers then any open cover of the unit space $G^0$ satisfies (P). If $G$ is ample and $\cal F = \{G^0\}$ then Example \ref{groupoid homology example} gives us an identification
\[ \bb L^{\cal F}_p \K_q(G \ltimes A) = \H_p(G;\K_q(A)). \]
\end{example}

\begin{rmk}\label{spatially weak equivalence remark}
Instead of considering the proper openly Morita embedded groupoids it may be desirable to consider a family $\cal F$ of open Morita embeddings of proper principal groupoids or of spaces. A condition on $\cal F$ analogous to (P) which considers open subsets of $G^0$ instead of arbitrary proper open subgroupoids of $G$ will ensure that $\cal F$ induces the same complementary pair as the family $\{G^0\}$. When $G$ has torsion in its isotropy this need no longer relate well to the topological K-theory $\Ktop_*(G)$ but for ample $G$ it will still be connected to the groupoid homology.
\end{rmk}

\begin{defn}
Let $\Omega \colon G \to H$ be an étale correspondence and let $\cal E$ and $\cal F$ be families of open Morita embeddings into $G$ and $H$ respectively. We say $\cal E$ and $\cal F$ are \textit{compatible} under $\Omega$ if for each embedding $K \hookrightarrow_M G$ in $\cal E$ there is an embedding $L \hookrightarrow_M H$ in $\cal F$ and an étale correspondence $\Lambda \colon K \to L$ such that 
\[ \begin{tikzcd}
G \arrow[r, "\Omega"] & H \\
K \arrow[u, hook] \arrow[r, "\Lambda"] & L \arrow[u, hook]
\end{tikzcd} \]
commutes up to isomorphism.
\end{defn}

An immediate consequence of compatibility of families $\cal E$ and $\cal F$ under $\Omega \colon G \to H$ is that $\Ind_\Omega(\I_{\cal F}) \subseteq \I_{\cal E}$.

\begin{example}\label{inverse semigroup equivariant topological correspondence example}
Let $S$ be an inverse semigroup and let $Z \colon X \to Y$ be an $S$-equivariant topological correspondence \cite[Definition 3.8]{Miller23a}. Explicitly, this consists of spaces $X$, $Y$ and $Z$ with actions of $S$, an $S$-equivariant continuous map $\rho \colon Z \to X$ and an $S$-equivariant local homeomorphism $\sigma \colon Z \to Y$ such that $\rho^{-1}(\dom_X s) = \dom_Z s$ holds for each $s \in S$. The resulting étale homomorphism $S \ltimes Z \to S \ltimes Y$ and action correspondence $S \ltimes X \to S \ltimes Z$ induce an étale correspondence $\Omega_{Z,S} \colon S \ltimes X \to S \ltimes Y$. Consider the families $\cal F^X_{\text{fin}}$ and $\cal F^Y_{\text{fin}}$ in $S \ltimes X$ and $S \ltimes Y$ respectively which satisfy (P) as in \Cref{inverse semigroup P example}. Their members are indexed by idempotents $e \in E(S)^\times$ and finite inverse subgroups $F \leq S_e$. For each such $(e,F)$ the associated members are the transformation groupoids $F_\downarrow \ltimes X$ and $F_\downarrow \ltimes Y$ by the inverse semigroup $F_\downarrow = \{ f d \mid d \in E(S), f \in F \}$ of restrictions of $F$ in $S$. Compatibility with $\Omega_{S,Z}$ follows from the commutativity up to isomorphism of
\[ \begin{tikzcd}
S \ltimes X \arrow[r, "\Omega_{S,Z}"] & S \ltimes Y \\
F_\downarrow \ltimes X \arrow[u, hook] \arrow[r, "\Omega_{F_\downarrow, Z}"] & F_\downarrow \ltimes Y. \arrow[u, hook]
\end{tikzcd} \]
\end{example}

\begin{prop}\label{ample unit space compatibility}
Let $\Omega \colon G \to H$ be a second countable Hausdorff étale correspondence of ample groupoids. Then $\{G^0\}$ and $\{H^0\}$ are compatible under $\Omega$. 
\begin{proof}
The local homeomorphism $q \colon \Omega \to \Omega/H$ has a continuous section $\lambda \colon \Omega/H \to \Omega$ because $\Omega$ is totally disconnected, second countable and Hausdorff. We may then construct an étale correspondence $\Lambda \colon G^0 \to H^0$ with bispace $\Omega/H$ with anchor maps $\orho \colon \Omega/H \to G^0$ and $\sigma \circ \lambda \colon \Omega/H \to H^0$. This makes
\[ \begin{tikzcd}
G \arrow[r, "\Omega"] & H \\
G^0 \arrow[u, hook] \arrow[r, "\Lambda"] & H^0 \arrow[u, hook]
\end{tikzcd} \]
commute up to isomorphism by construction.
\end{proof}
\end{prop}

With open Morita embeddings rather than just open subgroupoids, we can arrange for a compatible pair of families with respect to any étale correspondence.

\begin{prop}\label{transfer of OMEs}
Let $\Omega \colon G \to H$ be an \'etale correspondence and let $K \hookrightarrow_M G$ be an open Morita embedding into $G$. Then there is an open Morita embedding $L \hookrightarrow_M H$ and an actor $K \acts L$ with $K \bcdot L^0 = L$ such that 
\[ \begin{tikzcd}
G \arrow[r, "\Omega"] & H \\
K \arrow[u, hook] \arrow[r] & L \arrow[u, hook]
\end{tikzcd} \]
commutes. If $K$ is proper then so is $L$.
\begin{proof}
Let $\Lambda \colon K \to H$ be the composition $K \hookrightarrow_M G \xrightarrow{\Omega} H$. By \Cref{decomposition of correspondences} we may decompose $\Lambda$ into an actor $K \acts M$ and an open Morita embedding $M \hookrightarrow_M H$. The orbit $L = K \bcdot M^0$ of $M^0$ under the action $K \acts M$ is an open subgroupoid of $M$. The \'etale correspondence $\Lambda \colon K \to H$ therefore factorises into an actor $K \acts L$ and an open Morita embedding $L \hookrightarrow_M H$. Now suppose $K$ is proper and let $\rho \colon L^0 \to K^0$ be the anchor map for the actor. Then for any compact $C \subseteq L^0$, the restriction $L \restriction_C = (K \bcdot L^0) \restriction_C = K \restriction_{\rho(C)} \bcdot C$ is compact.
\end{proof}
\end{prop}

We note that if everything we started with was second countable or Hausdorff then so are the groupoids and correspondences constructed in this proof.

\begin{cor}\label{compatible family existence}
Let $G$ and $H$ be second countable Hausdorff étale groupoids and let $\Omega \colon G \to H$ be a second countable Hausdorff étale correspondence. Then $\Ind_\Omega(\fk N_H) \subseteq \fk N_G$. Furthermore, for any (countable) family $\cal E$ of proper second countable Hausdorff groupoids openly Morita embedded into $G$ satisfying (P) there is a (countable) family $\cal F$ of proper second countable Hausdorff groupoids openly Morita embedded into $H$ satisfying (P) which is compatible with $\cal E$ under $\Omega$. In particular, $\Ind_\Omega(\I_{\cal F}) \subseteq \I_{\cal E}$.
\begin{proof}
We first note that a countable family $\cal E$ of proper open subgroupoids satisfying (P) always exists by \Cref{inverse semigroup P example}, so $\Ind_\Omega(\fk N_H) \subseteq \fk N_G$ follows from the final part of the statement, by way of \Cref{properly weak equivalence result}.

For each $K \in \cal E$ use \Cref{transfer of OMEs} to pick a proper openly Morita embedded groupoid $L_K \hookrightarrow_M H$ and an étale correspondence $K \to L_K$ such that 
\[ \begin{tikzcd}
G \arrow[r, "\Omega"] & H \\
K \arrow[u, hook] \arrow[r] & L_K \arrow[u, hook]
\end{tikzcd} \]
commutes. Let $\cal F_0$ be a countable family of proper groupoids openly Morita embedded into $H$ satisfying (P) and set $\cal F = \{ L_K \mid K \in \cal E \} \cup \cal F_0$. Then $\cal F$ is compatible with $\cal E$ under $\Omega$ by construction, thus $\Ind(\I_{\cal F}) \subseteq \I_{\cal E}$, and $\cal F$ satisfies (P) because $\cal F_0$ does. 
\end{proof}
\end{cor}

This tells us that for any étale correspondence $\Omega \colon G \to H$ we can always build a compatible family $\cal F$ of embeddings into $H$ from a family $\cal E$ of embeddings into $G$, and to ensure that $\cal F$ satisfies (P) we can always artificially enlarge the family. In practice we might start with a pair of families we know are compatible, and ask when $\cal F$ inherits condition (P) from $\cal E$.

\begin{prop}\label{automatic P}
Let $\Omega \colon G \to H$ be an \'etale correspondence with $G \acts \Omega$ free such that for each $x \in H^0$ and finite isotropy subgroup $\Gamma \subseteq H^x_x$ there is $\omega \in \Omega$ with $\sigma(\omega) = x$ and $\omega \bcdot \Gamma \subseteq G \bcdot \omega$. Let $\cal F_G$ and $\cal F_H$ be families of proper groupoids openly Morita embedded into $G$ and $H$ which are compatible under $\Omega$. If $\cal F_G$ satisfies (P) then $\cal F_H$ satisfies (P).
\begin{proof}
We will use \Cref{checkable general P} to translate (P) for $\cal F_G$ into (P) for $\cal F_H$. Let $x \in H^0$ and let $\Gamma \leq H^x_x$ be a finite subgroup. By assumption there is $\omega_1 \in \Omega$ with $\sigma(\omega_1) = x$ and $\omega_1 \bcdot \Gamma \subseteq G \bcdot \omega$. By freeness of $G \acts \Omega$ there is a finite subgroup $\Gamma_1 \leq G^{\rho(\omega_1)}_{\rho(\omega_1)}$ with $\omega_1 \bcdot \Gamma = \Gamma_1 \bcdot \omega_1$. Using condition (P) for $\cal F_H$ there is a member $K \in \cal F_H$ with embedding $\Omega_K \colon K \hookrightarrow_M H$ and an element $\omega_2 \in \Omega_K$ such that $\sigma(\omega_2) = \rho(\omega_1)$ and $\omega_2 \bcdot \Gamma_1 \subseteq K \bcdot \omega_2$. There is therefore a finite subgroup $\Gamma_2 \leq K^{\rho(\omega_2)}_{\rho(\omega_2)}$ with $\omega_2 \bcdot \Gamma_1 = \Gamma_2 \bcdot \omega_2$. Pick an open Morita embedding $\Omega_L \colon L \hookrightarrow_M G$ and an \'etale correspondence $\Lambda \colon K \to L$ such that $\Omega \circ \Omega_K \cong \Omega_L \circ \Lambda$ and let $\Phi \colon \Omega_K \times_G \Omega \to \Lambda \times_L \Omega_L$ be the isomorphism of $K$-$H$-bispaces. Choose $\lambda \in \Lambda$ and $\omega_3 \in \Omega_L$ such that $[\lambda,\omega_3]_L = \Phi([\omega_2,\omega_1]_G)$. By construction, $[\Gamma_2 \bcdot \lambda, \omega_3]_L = [\lambda, \omega_3 \bcdot \Gamma]_L$. It follows that for each $h \in \Gamma$ there is $l \in L$ such that $\omega_3 \bcdot h = l \bcdot \omega_3$, i.e. $\omega_3 \bcdot \Gamma \subseteq L \bcdot \omega_3$.
\end{proof}
\end{prop}

\section{The interaction between Kasparov theory and groupoid homology}\label{groupoid homology section}

In \Cref{groupoid homology example} we see that for an ample groupoid $G$ the derived functors of crossed product K-theory functors on $\KK^G$ with respect to $\I_0 = \ker \Res^{G^0}_G$ are given by groupoid homology:
\[\bb L_p \K_q(G \ltimes -) \cong \H_p(G;\K_q(-)) \]
The goal of this section is to prepare ourselves to relate homological algebra in $\KK^G$ with respect to $\I_0$ to more classical homological algebra over $G$ through an analysis of the relationship between $G$-$\cs$-algebras $A \in \KK^G$ and the $G$-modules $\K_*(A) \in \lMod G$. For other ABC tuples $(\T,\I,F,A)$ it would be desirable to carry out a similar analysis to enable computation and handling of the derived functors $\bb L_n F (A)$.

\begin{defn}[Module over an ample groupoid]
Let $G$ be an ample groupoid. The \textit{groupoid ring} $\Z[G]$ is the abelian group of functions $G \to \Z$ generated by the indicator functions $\chi_U$ of compact open Hausdorff subsets $U \subseteq G$, with convolution as the ring\footnote{When $G^0$ is not compact $\Z[G]$ is not unital, but it is always locally unital.} multiplication. A (left) \textit{$G$-module} $M$ is a (left) $\Z[G]$-module such that $\Z[G] \bcdot M = M$. 
\end{defn}
For any totally disconnected locally LCH space $X$ we may consider the abelian group $\Z[X]$ as above, and an action $G \acts X$ induces a (left) $G$-module structure on $\Z[X]$. For an introduction to the category $\lMod G$ of $G$-modules we refer the reader to \cite[Section 2]{Miller24}. For an ample groupoid $G$ and a $G$-$\cs$-algebra $A = \Gamma_0(G^0, \A)$, the K-theory groups $\K_*(A)$ have the following $G$-module structure. For each clopen subset $V \subseteq G^0$ consider the subalgebra $VA \defeq C_0(V) A = \Gamma_0(V, \A)$. Inclusion and restriction maps $VA \to A \to VA$ induce injective and surjective maps $\K_*(VA) \to \K_*(A) \to \K_*(VA)$. By continuity of K-theory, $\K_*(A)$ is generated by the subgroups $\K_*(VA)$ for $V$ in any open cover of $G^0$. Let $\alpha \colon G \acts A$ denote the action of $G$ on $A$. Each compact slice $U \subseteq G$ induces a $*$-homomorphism $\alpha_U \colon s(U)A \to r(U)A$ defined at $a \in s(U)A$ by 
\[ \alpha_U(a) \colon x \mapsto \sum_{g \in G_x} \chi_U(g) (g \bcdot a_x). \]
Through restriction and inclusion this extends to a $*$-homomorphism $A \to A$ which we also call $\alpha_U$. It is straightforward to check by \cite[Lemma 2.4]{Miller24} that this induces the structure of a $G$-module on $K_i(A)$ for $i = 0,1$ with 
\[ \chi_U \cdot m = \K_i(\alpha_U)(m) \]
for each compact slice $U \subseteq G$ and $m \in \K_i(A)$. For second countable Hausdorff $G$, separable $G$-$\cs$-algebras $A$ and $B$ and $f \in \KK^G(A,B)$, the induced map $\K_i(f) \colon \K_i(A) \to \K_i(B)$ is $G$-equivariant, so we obtain a pair of functors 
\[\K_0, \K_1 \colon \KK^G \to \lMod G.\] 
The KK-theoretic induction functor of an étale correspondence has a module-theoretic analogue constructed by the author in \cite[Section 3]{Miller24}. For an étale correspondence $\Omega \colon G \to H$ of ample groupoids the induction functor 
\[\Ind_\Omega \colon \lMod H \to \lMod G\]
is given by the tensor product $\Z[\Omega] \otimes_H -$ over $\Z[H]$ with the $G$-$H$-bimodule $\Z[\Omega]$. To prepare ourselves to compare the two induction functors, suppose $G$, $H$ and $\Omega$ are Hausdorff, let $B$ be an $H$-$\cs$-algebra and let $U \subseteq \Omega$ be a compact slice. Let $\iota_U \colon \sigma(U)B \to B$ denote the inclusion. There is also an inclusion $\epsilon_U \colon \sigma(U) B \to \Ind_\Omega B$ given at $b \in \sigma(U) B$ by the element $\epsilon_U(b) \in \Ind_\Omega B \subseteq \Gamma_b(\Omega, \sigma^* \cal B)$ as follows for $b \in \sigma(U)B$ and $\omega \in \Omega$:
\begin{align*}
\epsilon_U \colon \sigma(U) B & \to \Ind_\Omega B \\
b & \mapsto \epsilon_U(b) \\
& \quad \; \;  \omega \mapsto \sum_{h \in H^{\sigma(\omega)}} \chi_U(\omega \cdot h) ( h \bcdot b_{s(h)} )
\end{align*}
This is the equivariant extension (see \cite[Proposition 4.3]{Miller23a}) of the section in $\Gamma_0(U, \sigma^* \B)$ determined by $b$ and the homeomorphism $\sigma \colon U \to \sigma(U)$. The image of $\epsilon_U$ is the subalgebra of functions supported on $U \cdot H \subseteq \Omega$, which is $q(U) \Ind_\Omega B$ with respect to the $C_0(\Omega/H)$-algebra structure on $\Ind_\Omega B$. We obtain inclusions 
\begin{align*}
 \K_*(\iota_U) \colon \K_*(\sigma(U)B) & \to \K_*(B) \\
\K_*(\epsilon_U) \colon \K_*(\sigma(U)B) & \to \K_*(\Ind_\Omega B) 
\end{align*}
whose images generate $\K_*(B)$ and $\K_*(\Ind_\Omega B)$ across all compact slices $U \subseteq \Omega$.

\begin{prop}[K-theory intertwines the induction functors]\label{K-theory intertwines}
Let $\Omega \colon G \to H$ be a second countable Hausdorff \'etale correspondence of second countable Hausdorff ample groupoids. Then for $B \in \KK^H$ there is a natural isomorphism of $G$-modules
\[\zeta_{\Omega,B} \colon \Ind_\Omega \K_*(B) \to \K_*(\Ind_\Omega B) \]
such that $\zeta_{\Omega,B}(\chi_U \otimes \K_*(\iota_U)(y)) = \K_*(\epsilon_U)(y)$ for each compact slice $U \subseteq \Omega$ and $y \in \K_*(\sigma(U)B)$.

\begin{proof}
Let $\alpha \colon G \acts \Ind_\Omega B$ and $\beta \colon H \acts B$ denote the respective groupoid actions on the $\cs$-algebras. For each compact slice $U \subseteq \Omega$ let $\eta_U \colon B \to \sigma(U) B$ denote the restriction. To construct $\zeta_{\Omega,B}$ we first construct a bilinear map
\[ L \colon \Z[\Omega] \times \K_*(B)  \to \K_*(\Ind_\Omega B) \]
such that $L(\chi_U, x) = \K_*(\epsilon_U \circ \eta_U)(x)$ for each compact slice $U \subseteq \Omega$ and $x \in \K_*(B)$. The assignment $(U,x) \mapsto \K_*(\epsilon_U \circ \eta_U)(x)$ respects disjoint unions of compact slices, so extends to a bilinear map $L$ by \cite[Lemma 2.4]{Miller24}. Let $V \subseteq H$ be a compact slice. For each compact slice $U \subseteq \Omega$ the diagram
\[ \begin{tikzcd}[column sep = large]
B \arrow[d, "\beta_V"] \arrow[rd, "\epsilon_{U \bcdot V} \circ \eta_{U \bcdot V}"]  \\
B \arrow[r, "\epsilon_U \circ \eta_U"'] & \Ind_\Omega B
\end{tikzcd} \]
commutes. It follows that $L$ is balanced and so extends to a homomorphism $\zeta_{\Omega,B} \colon \Ind_\Omega \K_*(B) \to \K_*(\Ind_\Omega B)$ of abelian groups. Moreover, the map $\zeta_{\Omega,B}$ is $G$-equivariant because for each pair of compact slices $U \subseteq \Omega$ and $W \subseteq G$ the diagram 
\[ \begin{tikzcd}[column sep = large]
B \arrow[r, "\epsilon_U \circ \eta_U"] \arrow[rd, "\epsilon_{W \bcdot U} \circ \eta_{W \bcdot U}"'] & \Ind_\Omega B \arrow[d, "\alpha_W"]  \\
   & \Ind_\Omega B 
\end{tikzcd} \]
commutes. Note that by the same argument $\zeta_{\Omega,B}$ is $\Omega/H$-equivariant, and thus it suffices to check bijectivity on the restrictions to $q(U) \subseteq \Omega/H$ for each compact slice $U \subseteq \Omega$. These restrictions are given by
\begin{align*}
\chi_{q(U)} \Ind_\Omega \K_*(B) & = \{ \chi_U \otimes \K_*(\iota_U)(y) \suchthat y \in \K_*(\sigma(U) B) \} \\
\chi_{q(U)} \K_*(\Ind_\Omega B) & = \{ \K_*(\epsilon_U)(y) \suchthat y \in \K_*(\sigma(U) B) \},
\end{align*}
upon which $\zeta_{\Omega,B}$ is bijective by construction.

Let $C \in \KK^H$ and consider a morphism $f \in \KK^H(B,C)$. For each compact slice $U \subseteq \Omega$ we may restrict $f$ and $\Ind_\Omega f$ to morphisms $\Res_{\sigma(U)} f \in \KK(\sigma(U) B, \sigma(U) C)$ and $\Res_{q(U)} \Ind_\Omega f \in \KK(q(U)\Ind_\Omega B, q(U)\Ind_\Omega C)$. To show naturality of $\zeta_{\Omega,-}$ we show commutativity of the following diagram:
\begin{equation}\label{diagram on naturality of intertwining}
\begin{tikzcd}
	{\sigma(U) B} && {q(U) \Ind_\Omega B} \\
	{\sigma(U) C} && {q(U) \Ind_\Omega C}
	\arrow["{\KK(\epsilon_U)}", "\cong"', from=1-1, to=1-3]
	\arrow["{\Res_{q(U)} \Ind_\Omega f}", from=1-3, to=2-3]
	\arrow["{\Res_{\sigma(U)} f}"', from=1-1, to=2-1]
	\arrow["{\KK(\epsilon_U)}", "\cong"', from=2-1, to=2-3]
\end{tikzcd}
\end{equation}
Let $(E,T)$ be an $H$-equivariant Kasparov cycle representing $f$. Consider the étale correspondence $U \bcdot H \colon q(U) \to H$. We may identify $q(U) \Ind_\Omega B$ with $\Ind_{U \bcdot H} B \restriction_{\sigma(U)}$ through the inclusion $\Gamma_b(U \cdot H, \sigma^* \B) \subseteq \Gamma_b(\Omega, \sigma^* \B)$. Similarly, $\Res_{q(U)} \Ind_\Omega f$ may be identified with $\Ind_{U \bcdot H} f$ (after the forgetful functor $\KK^{q(U)} \to \KK$) since each cutoff function for $\Omega$ restricts to a cutoff function for $U \bcdot H$. The correspondence $U \bcdot H \colon q(U) \to H \restriction_{\sigma(U)}$ has a cutoff function given by $\chi_U \colon U \bcdot H \to \bb R$, and so 
\[ \Ind_{U \bcdot H} f = [\Ind_{U \bcdot H} E, \Ind_{U \bcdot H, \chi_U} T]. \]
The $*$-isomorphisms labelled $\epsilon_U$ in (\ref{diagram on naturality of intertwining}) are compatible with the analogous isomorphism $\epsilon_U \colon \sigma(U) E \to \Ind_{U \bcdot H} E$, through which we may identify the correspondence $\sigma(U) E \colon \sigma(U) B \to \sigma(U) C$ with $\Ind_{U \bcdot H} E \colon \Ind_{U \bcdot H} B \to \Ind_{U \bcdot H} C$. This isomorphism conjugates the Fredholm operator $\sigma(U) T$ to $\Ind_{U \bcdot H, \chi_U} T$, from which it follows that (\ref{diagram on naturality of intertwining}) commutes.
\end{proof}
\end{prop}

For an ample groupoid $G$ and a $G$-module $M$ the \textit{coinvariants} $M_G$ is the tensor product $\Z[G^0] \otimes_G M$ over $\Z[G]$. This may be considered as the quotient of $M$ by the relations $\chi_U \bcdot m \sim \chi_{s(U)} \bcdot m$ for $m \in M$ and compact slices $U \subseteq G$. The quotient map $\pi_G \colon M \to M_G$ satisfies $\pi_G(m) = \chi_V \otimes m$ for each $m \in M$ and each compact open subset $V \subseteq G^0$ with $\chi_V \bcdot m = m$. Now consider an action $\alpha \colon G \acts A$ on a $\cs$-algebra $A$ and the inclusion map $\eta_A \colon A \to G \ltimes A$. Given an element $x \in \K_*(A)$ and a compact slice $U \subseteq G$, we have an equality $\K_*(\eta_A)(\chi_U \bcdot x) = \K_*(\eta_A)(\chi_{s(U)} \bcdot x)$. This is because the partial isometry $\chi_U \in M(G \ltimes A)$ conjugates $\eta_A \circ \alpha_U \colon A \to G \ltimes A$ to $\eta_A \circ \alpha_{s(U)} \colon A \to G \ltimes A$ in the sense that for each $a \in A$ we have
\begin{align*}
\chi_U^* \eta_A( \alpha_U(a)) \chi_U & = \eta_A(\alpha_{s(U)}(a)), \\
\chi_U \eta_A( \alpha_{s(U)}(a)) \chi_U^* & = \eta_A(\alpha_{U}(a)).
\end{align*}
It follows that $\K_*(\eta_A) \colon \K_*(A) \to \K_*(G \ltimes A)$ factors through the quotient map $\pi_G \colon \K_*(A) \to \K_*(A)_G$.
\begin{defn}[The comparison map]
Let $G$ be an ample groupoid and let $A$ be a $G$-$\cs$-algebra. The \textit{comparison map} $\gamma_A \colon \K_*(A)_G \to \K_*(G \ltimes A)$ is the homomorphism induced by the inclusion $A \subseteq G \ltimes A$. 
\end{defn}
For second countable Hausdorff $G$ the comparison map $\gamma_A$ is natural in $A \in \KK^G$ by naturality of the inclusion $A \subseteq G \ltimes A$ at the level of $\KK$. 

For an étale correspondence $\Omega \colon G \to H$ the author introduced \cite[Definition 6.10]{Miller23a} a natural transformation 
\[ \alpha_\Omega \colon \K_*(G \ltimes \Ind_\Omega -) \Rightarrow \K_*(H \ltimes -) \]
which relates the crossed product K-theory functors on $\KK^G$ and $\KK^H$. For modules of ample groupoids, the author also constructed \cite[Proposition 3.3]{Miller24} an $H$-equivariant homomorphism $\delta_\Omega \colon \Z[G^0] \otimes_G \Z[\Omega] \to \Z[H^0]$ which induces a natural transformation 
\[ \delta_\Omega \otimes - \colon (\Ind_\Omega -)_G \Rightarrow (-)_H \]
relating the coinvariants on $\lMod G$ and $\lMod H$.

\begin{prop}\label{modules and crossed product comparison}
Let $G$ and $H$ be second countable Hausdorff ample groupoids, let $\Omega \colon G \to H$ be a second countable Hausdorff étale correspondence and let $B \in \KK^H$ be an $H$-C*-algebra. Then the following diagram commutes:
\[\begin{tikzcd}
	{(\Ind_\Omega \K_*(B))_G} && {\K_*(B)_H} \\
	{\K_*(\Ind_\Omega B)_G} \\
	{\K_*(G \ltimes \Ind_\Omega B)} && {\K_*(H \ltimes B)}
	\arrow["{(\zeta_{\Omega,B})_G}", "\cong"', from=1-1, to=2-1]
	\arrow["{\gamma_{\Ind_\Omega B}}", from=2-1, to=3-1]
	\arrow["{\alpha_\Omega(B)}", from=3-1, to=3-3]
	\arrow["{\gamma_B}"', from=1-3, to=3-3]
	\arrow["{\delta_\Omega \otimes \id}"', from=1-1, to=1-3]
\end{tikzcd}\]
\begin{proof}
Fix a compact slice $U \subseteq \Omega$ and an element $y \in \K_*(\sigma(U) B)$. We need only check commutativity at the element $[\chi_U \otimes \K_*(\iota_U)(y)]_G \in (\Ind_\Omega \K_*(B))_G$. The homomorphism $\alpha_\Omega(B)$ is induced by a proper $\cs$-correspondence 
\[\Omega \ltimes B \colon G \ltimes \Ind_\Omega B \to H \ltimes B\]
introduced in \cite[Section 5]{Miller23a}. Consider this within the following diagram of proper $\cs$-correspondences with all other maps induced by $*$-homomorphisms:
\[\begin{tikzcd}
	{\sigma(U) B} && B \\
	{\Ind_\Omega B} \\
	{G \ltimes \Ind_\Omega B} && {H \ltimes B}
	\arrow["{\epsilon_U}", from=1-1, to=2-1]
	\arrow[hook, from=2-1, to=3-1]
	\arrow["{\Omega \ltimes B}", from=3-1, to=3-3]
	\arrow[hook, from=1-3, to=3-3]
	\arrow[hook, "\iota_U", from=1-1, to=1-3]
\end{tikzcd}\]
By construction, the resulting diagram in K-theory applied to $y$ computes the values of the initial diagram applied to $[\chi_U \otimes \K_*(\iota_U)(y)]_G$. It therefore suffices to show that the above diagram commutes in the category of $\cs$-correspondences. The $\downthenrightarrow$ composition of $\cs$-correspondences is given by the Hilbert $H \ltimes B$-submodule $\overline{\Gamma_c(U \bcdot H, \sigma^* \B)}$ of $\Omega \ltimes B = \overline{\Gamma_c(\Omega,\sigma^* \B)}$ with the action of $b \in \sigma(U)B$ as follows on $\xi \in \Gamma_c(U \bcdot H, \sigma^* \B)$, for $\omega \in \Omega$:
\begin{align*}
\sigma(U) B \times \Gamma_c(U \bcdot H, \sigma^* \B) & \to \Gamma_c(U \bcdot H, \sigma^* \B) \\
(b,\xi) & \mapsto \epsilon_U(b) \bcdot \xi \\
& \quad \; \;  \omega \mapsto \sum_{h \in H^{\sigma(\omega)}} \chi_U(\omega \cdot h) ( h \bcdot b_{s(h)} ) \xi(\omega)
\end{align*}
On the other hand, the $\rightthendownarrow$ composition is simply the inclusion $\sigma(U) B \subseteq H \ltimes B$, which as a $\cs$-correspondence has underlying Hilbert $H \ltimes B$-module $ \overline{\Gamma_c(H^{\sigma(U)}, s^* \B)}$, as a submodule of $H \ltimes B = \overline{\Gamma_c(H,s^* \B)}$. The homeomorphism $U \bcdot H \to H^{\sigma(U)}$ satisfying $u \bcdot h \mapsto h$ for $u \in U$ and $h \in H$ with $\sigma(u) = r(h)$ induces an isomorphism of these $\cs$-correspondences.
\end{proof}
\end{prop}

\begin{prop}\label{comparison map isomorphism}
Let $G$ be a second countable Hausdorff ample groupoid, let $B \in \KK^{G^0}$ be a separable $C_0(G^0)$-algebra and consider the induced algebra $A = \Ind^G_{G^0} B \in \KK^G$. Then the comparison map $\gamma_A \colon \K_*(A)_G \to \K_*(G \ltimes A)$ is an isomorphism.
\begin{proof}
We may consider $A$ as the induced algebra $\Ind_\Omega B$ for the étale correspondence $\Omega = G \colon G \to G^0$. Applying \Cref{modules and crossed product comparison} we obtain the following commutative diagram.
\[\begin{tikzcd}
	{(\Ind_\Omega \K_*(B))_G} && {\K_*(B)} \\
	{\K_*(A)_G} \\
	{\K_*(G \ltimes A)} && {\K_*(B)}
	\arrow["{(\zeta_{\Omega,B})_G}", "\cong"', from=1-1, to=2-1]
	\arrow["{\gamma_{A}}", from=2-1, to=3-1]
	\arrow["{\alpha_{\Omega}(B)}", from=3-1, to=3-3]
	\arrow["{\id}"', from=1-3, to=3-3]
	\arrow["{\delta_\Omega \otimes \id}"', from=1-1, to=1-3]
\end{tikzcd}\]
The map $\delta_\Omega \colon \Z[G^0] \otimes_G \Z[G] \to \Z[G^0]$ is given by $\delta_\Omega(\eta \otimes \xi) = s_*( \eta \vartriangleright \xi) = \eta \vartriangleleft \xi$ for $\eta \in \Z[G^0]$ and $\xi \in \Z[G]$, where $\vartriangleright$ denotes the action of $\Z[G^0]$ on $\Z[G]$ and $\vartriangleleft$ vice versa. It follows that $\delta_\Omega$ and therefore $\delta_\Omega \otimes \id$ is an isomorphism. The map $\alpha_\Omega(B)$ is induced by the proper correspondence $\Omega \ltimes B \colon G \ltimes A \to B$. The action of $G$ on $A$ factors through an action of $G \ltimes G$, so we may identify $G \ltimes A$ with the crossed product $(G \ltimes G) \ltimes A$. Consider the invertible étale correspondence $\Lambda = G \colon G \ltimes G \to G^0$ and note that $\Ind_\Lambda B = A$. We may identify $\Omega \ltimes B \colon G \ltimes A \to B$ with $\Lambda \ltimes B \colon (G \ltimes G) \ltimes A \to B$, which inherits invertibility from $\Lambda$ by \cite[Proposition 7.4]{Miller23a}. Since all other maps in the diagram are isomorphisms, $\gamma_A$ must also be.
\end{proof}
\end{prop}

\begin{rmk}\label{homology identified as a derived functor}
Let $G$ be a second countable Hausdorff ample groupoid, let $A \in \KK^G$ be a $G$-$\cs$-algebra and consider the ideal $\I_0$ as in \Cref{groupoid homology example}. Then the comparison map for induced algebras induces an isomorphism
\[\H_p(G;\K_q(A)) \cong \bb L_p \K_q(G \ltimes A)\]
as follows. For $n \geq 0$ let $P_n = (\Ind^G_{G^0} \Res^{G^0}_G)^{n+1} A$ and consider the $\I_0$-projective resolution $P_\bullet \to A$ induced by the induction-restriction adjunction $\Ind^G_{G^0} \dashv \Res^{G^0}_G$ as in \Cref{adjunctions are good}. The homology groups for the chain complex $\K_q(G \ltimes P_\bullet)$ are $\bb L_* \K_q(G \ltimes A)$ as $P_\bullet \to A$ is $\I_0$-projective. Furthermore, by \Cref{K-theory intertwines} the $G$-modules $\K_q(P_n)$ are isomorphic to $G^0$-induced modules, and so by \cite[Corollary 2.20]{Miller24} are $G$-acyclic. The homology of the chain complex $\K_q(P_\bullet)_G$ therefore produces $\H_*(G; \K_q(A))$. By naturality the comparison maps
\[ \gamma_{P_\bullet} \colon \K_q(P_\bullet)_G \to \K_q(G \ltimes P_\bullet) \]
form a chain map. Each $\gamma_{P_n}$ is an isomorphism by \Cref{comparison map isomorphism}, so they in particular induce an isomorphism of the homology groups.
\end{rmk}

\section{The ABC category}\label{ABC section}

We are now ready to set up a category of ABC tuples $(\T,\I,F,A)$ so that we can study functoriality properties of the localisation $\bb L^\I F(A)$, the derived functors $\bb L^\I_n F(A)$ and ultimately the associated ABC spectral sequence. Our starting point is to consider a graded homomorphism $\varphi \colon F(A) \to F'(A')$ relating to a second ABC tuple $(\T',\I',F',A')$. Morphisms of ABC tuples should be triples that lift $\varphi$ to the level of triangulated categories, that is
\begin{itemize}
\item a triangulated functor $\Phi \colon \T' \to \T$,
\item a natural transformation $\alpha \colon F_* \circ \Phi \Rightarrow F'_*$,
\item and a morphism $f \colon A \to \Phi(A')$,
\end{itemize}
such that $\varphi = \alpha_{A'} \circ F_*(f)$. We write $\App(\Phi, \alpha, f)$ for the composition $\alpha_{A'} \circ F_*(f)$.
\[\begin{tikzcd}
	{F(A)} && {F'(A')} \\
	& {F(\Phi (A'))}
	\arrow["{\alpha_{A'}}"', from=2-2, to=1-3]
	\arrow["{F(f)}"', from=1-1, to=2-2]
	\arrow["{\App(\Phi,\alpha,f)}", from=1-1, to=1-3]
\end{tikzcd}\]

\begin{example}\label{basic correspondence example}
Let $G$ and $H$ be second countable Hausdorff \'etale groupoids and let $\Omega \colon G \to H$ be a proper second countable Hausdorff \'etale correspondence. The $\cs$-correspondence $\cs(\Omega) \colon \cs(G) \to \cs(H)$ induces a map in K-theory
\[ \K_*(\cs(\Omega)) \colon \K_*(\cs(G)) \to \K_*(\cs(H)). \]
Consider the categories $\KK^G$ and $\KK^H$, the functors $\K_*(G \ltimes -) \colon \KK^G \to \Ab_*$ and $\K_*(H \ltimes -) \colon \KK^H \to \Ab_*$ and the objects $C_0(G^0)$ and $C_0(H^0)$. In \cite{Miller23a} the author constructed a triple consisting of
\begin{itemize}
\item the induction functor $\Ind_\Omega \colon \KK^H \to \KK^G$,
\item the induction natural transformation $\alpha_\Omega \colon \K_*(G \ltimes \Ind_\Omega -) \Rightarrow \K_*(H \ltimes -)$,
\item and the morphism $f_\Omega \colon C_0(G^0) \to \Ind_\Omega C_0(H^0)$ in $\KK^G$ induced by the $G$-equivariant proper map $\orho \colon \Omega/H \to G^0$,
\end{itemize}
such that $\App(\Ind_\Omega, \alpha_\Omega, f_\Omega) = \K_*(\cs(\Omega)) $. More generally we can drop the properness requirement on $\Omega$ and consider coefficient $\cs$-algebras $A \in \KK^G$ and $B \in \KK^H$ so long as we have a morphism $f \in \KK^G(A,\Ind_\Omega B)$. We then set 
\[\K_*(\Omega;f) \defeq \App(\Ind_\Omega, \alpha_\Omega, f) \colon \K_*(G \ltimes A) \to \K_*(H \ltimes B).\]
\end{example}
Triples of this form make sense and appear outside of the triangulated setting too:
\begin{example}\label{basic homology correspondence example}
Let $G$ and $H$ be ample groupoids and let $\Omega \colon G \to H$ be a proper \'etale correspondence. In \cite{Miller24} the author constructed a graded homomorphism
\[ \H_*(\Omega) \colon \H_*(G) \to \H_*(H) \]
as the application $\App(\Ind_\Omega, \delta_\Omega \otimes \id, \orho^*)$ of a triple on the level of the module categories $\lMod G$ and $\lMod H$: 
\begin{itemize}
\item the induction functor $\Ind_\Omega \colon \lMod H \to \lMod G$,
\item the natural transformation $\delta_\Omega \otimes \id \colon (\Ind_\Omega -)_G \Rightarrow (-)_H$,
\item and the morphism $\orho^* \colon \Z[G^0] \to \Ind_\Omega \Z[H^0]$  induced by the $G$-equivariant proper map $\orho \colon \Omega/H \to G^0$.
\end{itemize}

\end{example}
\begin{rmk}
It is possible to place the classical homological algebra of abelian categories such as $\lMod G$ within the framework of triangulated categories with homological ideals by considering the (triangulated) homotopy category $\Ho(\lMod G)$ with the kernel of the homology functor $\H_* \colon \Ho(\lMod G) \to (\lMod G)_*$. The above example can then be extended to a triple at the triangulated categories level. Furthermore, the K-theory groups induce an additive (but not triangulated!) functor $\KK^G \to \Ho(\lMod G)$ which relates the above two examples.
\end{rmk}

\begin{defn}[Composition of triples] 
Let $\T$, $\T'$ and $\T''$ be categories with functors $F$, $F'$ and $F''$ to a common target category $\cal C$ and objects $A$, $A'$ and $A''$. Let $\Phi \colon \T' \to \T$ and $\Psi \colon \T'' \to \T'$ be functors, with natural transformations $\alpha \colon F \circ \Phi \Rightarrow F'$ and $\beta \colon F' \circ \Psi \Rightarrow F''$ and morphisms $f \colon A \to \Phi A'$ and $g \colon A' \to \Psi A''$. The \textit{composition} $(\Psi, \beta, g) \circ (\Phi, \alpha, f)$ consists of:
\begin{itemize}
\item the functor $\Phi \circ \Psi \colon  \T'' \rightarrow \T$,
\item the natural transformation $\beta \circ (\alpha \Psi) \colon F \circ \Phi \circ \Psi \Rightarrow F''$, 
\item and the morphism $\Phi(g) \circ f \colon A \to \Phi \Psi A''$. 
\end{itemize}
\end{defn}
\begin{rmk}
The assignment $(\Phi, \alpha, f) \mapsto \App(\Phi,\alpha,f) = \alpha_{A'} \circ F(f)$ is functorial with respect to this composition. 
\end{rmk}
Given proper \'etale correspondences $\Omega \colon G \to H$ and $\Lambda \colon H \to K$ as in \Cref{basic correspondence example}, the composition $(\Ind_\Lambda, \alpha_\Lambda, f_\Lambda) \circ (\Ind_\Omega, \alpha_\Omega, f_\Omega) $ is equivalent to the triple $(\Ind_{\Lambda \circ \Omega}, \alpha_{\Lambda \circ \Omega}, f_{\Lambda \circ \Omega})$ in the following sense. There is a natural isomorphism 
\[ \Ind_\Omega \circ \Ind_\Lambda \cong \Ind_{\Lambda \circ \Omega}\]
which identifies $\alpha_\Lambda \circ (\alpha_\Omega \Ind_\Lambda)$ with $\alpha_{\Lambda \circ \Omega}$ and $\Ind_\Omega(f_\Lambda) \circ f_\Omega$ with $f_{\Lambda \circ \Omega}$ \cite[Section 7]{Miller23a}.
\begin{defn}[Equivalence of triples]
Consider categories $\T$ and $\T'$ with functors $F$ and $F'$ to a common target category $\cal C$ and objects $A$ and $A'$. For $i = 1,2$ let $\Phi_i \colon \T' \to \T$ be a functor, $\alpha_i \colon F \circ \Phi_i \Rightarrow F'$ a natural transformation and $f_i \colon A \to \Phi_i A'$ a morphism. We say the triples $(\Phi_1,\alpha_1,f_1)$ and $(\Phi_2, \alpha_2, f_2)$ are \textit{equivalent} if there is a natural isomorphism $\eta \colon \Phi_1 \cong \Phi_2$ such that $\alpha_1 = \alpha_2 \circ F \eta$ and $\eta_{A'} f_1 = f_2$.
\end{defn}

When we work with triangulated categories $\T$ and $\T'$ with distinguished homological ideals $\I \vartriangleleft \T$ and $\I' \vartriangleleft \T'$, we typically want a functor $\Phi \colon \T' \to \T$ to be triangulated and map $\I'$ into $\I$. It then automatically sends $\I'$-exact sequences to $\I$-exact sequences and $\I'$-equivalences to $\I$-equivalences. For an étale correspondence $\Omega \colon G \to H$ with ideals $\I_{\cal E} \vartriangleleft \KK^G$ and $\I_{\cal F} \vartriangleleft \KK^H$ defined by families $\cal E$ and $\cal F$ of open Morita embeddings into $G$ and $H$, $\Ind_\Omega(\I_{\cal F}) \subseteq \I_{\cal E}$ if $\cal E$ and $\cal F$ are compatible under $\Omega$. For details on the construction of compatible families see \Cref{inverse semigroup equivariant topological correspondence example}, \Cref{ample unit space compatibility} and \Cref{compatible family existence}. 

\begin{defn}[Morphism of ABC tuples]
Given ABC tuples $ \fk{M} =(\T,\I,F,A)$ and $ \fk{M}' = (\T',\I',F',A')$, an \textit{ABC cycle} from $\fk M$ to $\fk M '$ is given by a triple $(\Phi,  \alpha, f)$:
\begin{itemize}
\item a triangulated functor $\Phi \colon  \T' \rightarrow \T$ mapping $\I'$ into $\I$,
\item a natural transformation $\alpha \colon  F \circ \Phi \Rightarrow F'$,
\item and a morphism $f \colon  A \rightarrow \Phi(A')$ in $\T$.
\end{itemize}
An \textit{ABC morphism} $\fk m \colon \fk M \to \fk M'$ is an equivalence class $\fk m = [\Phi, \alpha, f]$ of ABC cycles under equivalence of triples. We will usually work directly with ABC cycles with the understanding that our constructions depend only on the equivalence class.
\end{defn}

\begin{example}\label{motivating ABC morphism example}
Let $G$ and $H$ be second countable Hausdorff étale groupoids with coefficient $\cs$-algebras $A \in \KK^G$ and $B \in \KK^H$ and countable families $\cal E$ and $\cal F$ of open Morita embeddings into $G$ and $H$ respectively. Let $\Omega \colon G \to H$ be a second countable Hausdorff étale correspondence under which $\cal E$ and $\cal F$ are compatible and let $f \in \KK^G(A,\Ind_\Omega B)$ be a morphism. Then 
\[ [\Ind_\Omega, \alpha_\Omega, f] \colon (\KK^G,\I_{\cal E}, \K_*(G \ltimes -), A) \to (\KK^H,\I_{\cal F}, \K_*(H \ltimes -), B) \]
is an ABC morphism. This has application $\K_*(\Omega;f) \colon \K_*(G \ltimes A) \to \K_*(H \ltimes B)$.
\end{example}

We now consider the functoriality of the localisation functors. Recall that localisation makes sense for any complementary pair $(\fk P, \fk N)$ of subcategories of a triangulated category $\T$. A morphism is a \textit{weak equivalence} if its cone is in $\fk N$.

\begin{lem}\label{localisation morphism lemma}
Let $\T$ and $\T'$ be triangulated categories with complementary pairs of subcategories with associated localisation functors $L \colon \T \to \T$ and $L' \colon \T' \to \T'$ and localisation natural transformations $\mu \colon L \Rightarrow \id_\T$ and $\mu' \colon L' \Rightarrow \id_{\T'}$. Let $\mu \colon L \Rightarrow \id_\T$ and $\mu' \colon L' \Rightarrow \id_{\T'}$ denote the localisation natural transformations and let $\Phi \colon \T' \to \T$ be a functor sending weak equivalences to weak equivalences. Then there exists a unique natural transformation $\bb L_\Phi \colon L \Phi \Rightarrow \Phi L'$ such that the following diagram commutes.
\begin{equation}\label{localisaion lemma diagram}
\begin{tikzcd}
	{L \Phi } && {\Phi L'} \\
	\\
	 && \Phi
	\arrow["\Phi \mu'", Rightarrow, from=1-3, to=3-3]
	\arrow["\mu \Phi"', Rightarrow, from=1-1, to=3-3]
	\arrow[Rightarrow, from=1-1, to=1-3, "\bb L_\Phi"]
\end{tikzcd}
\end{equation}
For $\Phi = \id_\T$ and $L' = L$, we have $\bb L_{\id_\T} = \id_L$. If $\T''$ is another triangulated category with a pair of complementary subcategories and a functor $\Psi \colon \T'' \to \T'$ maps weak equivalences to weak equivalences, then $\bb L_{\Phi \circ \Psi} =  (\Phi \bb L_\Psi) \circ (\bb L_\Phi \Psi)$. 
\begin{proof}
The localisation functor $L$ maps weak equivalences to isomorphisms, so $L \Phi \mu'$ is invertible. Applying the naturality of $\mu$ to $\Phi \mu'$ this fits into a commutative square:
\[\begin{tikzcd}
	{L \Phi L' } && {L \Phi } \\
	\\
	{\Phi L'} && \Phi
	\arrow["\Phi \mu'", Rightarrow, from=3-1, to=3-3]
	\arrow["\mu \Phi"', Rightarrow, from=1-3, to=3-3]
	\arrow["\mu \Phi L'"', Rightarrow, from=1-1, to=3-1]
	\arrow["L \Phi \mu'", Rightarrow, from=1-1, to=1-3]
	\arrow[red, Rightarrow, from=1-3, to=3-1, "\bb L_\Phi"]
\end{tikzcd}\]
We set $\bb L_\Phi = (\mu \Phi L') \circ (L \Phi \mu')^{-1}$. For uniqueness, note first that \eqref{localisaion lemma diagram} determines $\bb L_\Phi L'$ as $(\Phi \mu' L')^{-1} \circ \mu \Phi L'$. Applying naturality of $\bb L_\Phi$ to $\mu'$ we obtain the commutative square:
\[\begin{tikzcd}
	{L \Phi L' } && {L \Phi } \\
	\\
	{\Phi L' L'} && \Phi L'
	\arrow["\Phi L' \mu'", Rightarrow, from=3-1, to=3-3]
	\arrow["\bb L_\Phi"', Rightarrow, from=1-3, to=3-3]
	\arrow["\bb L_\Phi L'"', Rightarrow, from=1-1, to=3-1]
	\arrow["L \Phi \mu'", Rightarrow, from=1-1, to=1-3]
\end{tikzcd}\]
Invertibility of $L \Phi \mu'$ implies that $\bb L_\Phi$ is determined by $\bb L_\Phi L'$ and thus unique. Compatibility with identities and with composition follow straightforwardly from uniqueness.
\end{proof}
\end{lem}

\begin{defn}\label{localisation morphism definition}
Let $\T$ and $\T'$ be triangulated categories with homological ideals $\I$ and $\I'$, functors $F \colon \T \to \fk C$ and $F' \colon \T' \to \fk C$ and objects $A \in \T$ and $A' \in \T'$. Let $\Phi \colon \T' \to \T$ be a functor sending $\I'$-equivalences to $\I$-equivalences, $\alpha \colon F \circ \Phi \Rightarrow F'$ a natural transformation and $f \colon A \to \Phi(A')$ a morphism. The \textit{localisation application morphism} $\bb L(\Phi, \alpha, f)$ is the morphism
\[ \bb L(\Phi, \alpha, f) = \App(\Phi, \alpha, \bb L_\Phi(A') \circ L(f)) \colon \bb L^\I F(A) \to \bb L^{\I'} F'(A').  \]
\end{defn}
By \Cref{localisation morphism lemma} this is compatible with composition of triples and fits into the following diagram.
\[\begin{tikzcd}
	{\bb L^{\I} F(A)} && {\bb L^{\I'} F'(A')} \\
	{F(A)} && {F'(A')}
	\arrow["{\bb L(\Phi, \alpha, f)}", from=1-1, to=1-3]
	\arrow["{F'(\mu'_{A'})}", from=1-3, to=2-3]
	\arrow["{F(\mu_A)}"', from=1-1, to=2-1]
	\arrow["{\App(\Phi, \alpha, f)}", from=2-1, to=2-3]
\end{tikzcd}\]

\begin{example}\label{functoriality of topological K theory}
Let $G$ and $H$ be second countable Hausdorff étale groupoids with coefficient $\cs$-algebras $A \in \KK^G$ and $B \in \KK^H$, let $\Omega \colon G \to H$ be a second countable Hausdorff étale correspondence and let $f \in \KK^G(A,\Ind_\Omega B)$. Recall that with respect to the complementary pairs $(\fk P_G, \fk N_G)$ and $(\fk P_H, \fk N_H)$ the localisations are given by topological K-theory (\Cref{baum connes assembly map identification}, \Cref{full reduced remark}): 
\begin{align*}
\bb L \K_*(G \ltimes A) & \cong \Ktop_*(G;A) \\
\bb L \K_*(H \ltimes B) & \cong \Ktop_*(H;B)
\end{align*}
The localisation application morphism for the triple $(\Ind_\Omega, \alpha_\Omega, f)$ gives rise to functoriality of the topological K-theory:
\[ \Ktop_*(\Omega;f) \defeq \bb L(\Ind_\Omega, \alpha_\Omega, f) \colon \Ktop_*(G;A) \to \Ktop_*(H;B) \]
This is compatible with the full assembly maps. To relate it to the Baum--Connes assembly maps (i.e., the reduced assembly maps), first consider the proper $\cs$-correspondence $\Omega \ltimes B \colon G \ltimes \Ind_\Omega B \to H \ltimes B$ which induces $\alpha_\Omega$ (see \cite[Section 5]{Miller23a}). The quotient $H \ltimes B \to H \ltimes_r B$ induces a reduced Hilbert $H \ltimes_r B$-module $\Omega \ltimes_r B$. Suppose that the $*$-homomorphism $G \ltimes \Ind_\Omega B \to \cal K(\Omega \ltimes_r B)$ factors through the quotient map $G \ltimes \Ind_\Omega B \to G \ltimes_r \Ind_\Omega B$. Through the resulting proper $\cs$-correspondence $\Omega \ltimes_r B \colon G \ltimes_r \Ind_\Omega B \to H \ltimes_r B$ we may define a map $\Kred_*(\Omega;f) \defeq \K_*(\Omega \ltimes_r B) \circ \K_*(G \ltimes_r f)$ in K-theory at the reduced level making the following diagram commute:
\[\begin{tikzcd}
	{\Ktop_*(G;A)} & {\K_*(G \ltimes A)} & {\K_*(G \ltimes_r A)} \\
	{\Ktop_*(H;B)} & {\K_*(H \ltimes B)} & {\K_*(H \ltimes_r B)}
	\arrow[from=1-1, to=1-2]
	\arrow["{\mu_{G,A}}"', from=1-1, to=1-3, bend left]
	\arrow["{\K_*(\Omega;f)}", from=1-2, to=2-2]
	\arrow["{\Ktop_*(\Omega;f)}", from=1-1, to=2-1]
	\arrow[from=2-1, to=2-2]
	\arrow["{\mu_{H,B}}", from=2-1, to=2-3, bend right]
	\arrow["{\Kred_*(\Omega;f)}", from=1-3, to=2-3]
	\arrow[from=1-2, to=1-3]
	\arrow[from=2-2, to=2-3]
\end{tikzcd}\]
In particular, if $(G,A)$ and $(H,B)$ satisfy the Baum--Connes conjecture with coefficients then $\Ktop_*(\Omega;f)$ may be identified with $\Kred_*(\Omega;f)$.
\end{example}

\begin{prop}[Derived functor maps of an ABC morphism]\label{derivedmaps}
Let $\T$ and $\T'$ be triangulated categories with homological ideals $\I$ and $\I'$, stable homological functors $F$ and $F'$ into $\Ab_*$ and objects $A$ and $A'$. Let $\Phi \colon \T' \to \T$ be an additive functor sending $\I'$-exact sequences to $\I$-exact sequences, let $\alpha \colon F \circ \Phi \Rightarrow F'$ be a natural transformation and let $f \colon A \to F(A')$ be a morphism. Then for each $n \geq 0$ there is a homomorphism 
\[\bb L_n(\Phi, \alpha, f) \colon \bb L_n F(A) \to \bb L_n F'(A')\]
with the following property. Let $P_\bullet \to A$ and $P'_\bullet \to A'$ be projective resolutions with respect to $\I$ and $\I'$ respectively and let $ \tilde{f} \colon P_\bullet \to \Phi(P'_\bullet)$ be a chain map over $f \colon A \to \Phi(A')$. Let $\alpha_{P'} \colon F\Phi(P'_\bullet) \to F'(P'_\bullet)$ be the chain map induced by $\alpha$. Then the chain map $\alpha_{P'} \circ F (\tilde{f}) \colon F(P_\bullet) \to F'(P'_\bullet)$ over $\App(\Phi, \alpha, f) \colon F(A) \to F'(A')$ induces $\bb L_n(\Phi, \alpha, f) \colon \bb L_n F(A) \to \bb L_n F'(A')$ by taking homology.
\[ \begin{tikzcd}
\cdots \arrow[r] & F(P_1) \arrow[r] \arrow[d, "F(\tilde{f}_1)"]     & F(P_0) \arrow[r] \arrow[d, "F(\tilde{f}_0)"]     & F(A) \arrow[d, "F(f)"] \arrow[r]             & 0 \\
\cdots \arrow[r] & F\Phi(P'_1) \arrow[r] \arrow[d, "\alpha_{P'_1}"] & F\Phi(P'_0) \arrow[r] \arrow[d, "\alpha_{P'_0}"] & F\Phi(A') \arrow[d, "\alpha_{A'}"] \arrow[r] & 0 \\
\cdots \arrow[r] & F'(P'_1) \arrow[r]                               & F'(P'_0) \arrow[r]                               & F'(A') \arrow[r]                             & 0
\end{tikzcd} \]
Furthermore, the assignment $(\Phi, \alpha, f) \mapsto \bb L_n(\Phi, \alpha, f)$ respects identities, composition and equivalence of triples.

\begin{proof}
To construct $\bb L_n (\Phi, \alpha, f)$, we consider arbitrary projective resolutions $P_\bullet \to A$ and $P'_\bullet \to A'$ with respect to $\I$ and $\I'$ respectively. By \Cref{resolutionmap}, there is a lift $\tilde{f} \colon P_\bullet \to \Phi(P'_\bullet)$ of $f \colon A \to \Phi(A')$ which is unique up to chain homotopy. We may then define $\bb L_n (\Phi, \alpha, f)$ as the map in homology induced by the chain map $\alpha_{P'} \circ F(\tilde{f})$, which is a lift of $\App(\Phi, \alpha, f) = \alpha_{A'} \circ F(f)$.

To check that this is well-defined, suppose we have some other projective resolutions $Q_\bullet \to A$ and $Q'_\bullet \to A'$ and a lift $\hat{f} \colon Q_\bullet \to \Phi(Q'_\bullet)$ of $f$. Pick lifts $\tilde{\id} \colon P_\bullet \to Q_\bullet$ and $\tilde{\id}' \colon P'_\bullet \to Q'_\bullet$ of $\id \colon A \to A$ and $\id' \colon A' \to A'$. Then by \Cref{resolutionmap}, $\tilde{\id}' \circ \tilde{f}$ is chain homotopic to $\hat{f} \circ \tilde{\id}$, and therefore the following diagram commutes.
\[ \begin{tikzcd}[column sep = large]
H_n(F(P_\bullet)) \arrow[rr, "H_n\left(\alpha_{P'} \circ F(\tilde{f})\right)"] \arrow[d, "H_n(F(\tilde{\id}))", "\cong"'] &  & H_n(F'(P'_\bullet)) \arrow[d, "H_n(F'(\tilde{\id}'))", "\cong"'] \\
H_n(F(Q_\bullet)) \arrow[rr, "H_n\left(\alpha_{Q'} \circ F(\hat{f})\right)"]                                    &  & H_n(F'(Q'_\bullet))                                   
\end{tikzcd} \]
The compatibility of $(\Phi,\alpha,f) \mapsto \bb L_n(\Phi, \alpha, f)$ with identities, composition and equivalence each follow straightforwardly from well-definition.
\end{proof}
\end{prop}

\begin{rmk}\label{derived functor specialisation}
If we keep the triangulated category $\T$, homological ideal $\I$ and stable homological functor $F$ constant, an ABC morphism specialises to a morphism $f \colon A \to B$ in $\T$. It is straightforward to check that $\bb L(\id_\T, \id_F, f) = \bb L F(f) \colon \bb L F(A) \to \bb L F(B)$ and that $\bb L_n(\id_\T, \id_F, f) = \bb L_n F(f) \colon \bb L_n F(A) \to \bb L_n F(B)$ for each $n \geq 0$.
\end{rmk}

A more concrete description of $\bb L_n(\fk m)$ depends upon a more concrete description of the derived functors. In \Cref{groupoid homology section} we analysed the description of derived functors as groupoid homology in the setting of an ample groupoid (see in particular \Cref{homology identified as a derived functor}). The work in that section on compatibility with an étale correspondence $\Omega \colon G \to H$ allows us to describe the derived functor maps $\bb L_n(\fk m)$ as the induced map in homology $\H_n(\Omega;-)$ constructed by the author in \cite[Theorem 3.5]{Miller24}:
\begin{thm}\label{induced maps in homology agree}
Let $\Omega \colon G \to H$ be a second countable Hausdorff \'etale correspondence of second countable Hausdorff ample groupoids, let $A \in \KK^G$ and $B \in \KK^H$ and let $f \in \KK^G(A, \Ind_\Omega B)$. Consider the homological ideal $\I_0^G = \ker \Res^{G^0}_G \vartriangleleft \KK^G$ and the associated derived functors $\bb L_n$ for $n \geq 0$, and similarly for $H$. Recall the identification $\H_n(G ; \K_i(-)) \cong \bb L_n \K_i(G \ltimes -)$ in \Cref{homology identified as a derived functor}. For $i = 0,1$ let $\K_i(f) \colon \K_i(A) \to \Ind_\Omega \K_i(B)$ denote the induced map in K-theory with the identification $\K_*(\Ind_\Omega B) \cong \Ind_\Omega \K_*(B)$. Then the following diagram commutes:
\[\begin{tikzcd}[column sep = 3cm]
 \H_n(G; \K_i(A)) \arrow[r, "{\H_n(\Omega ; \K_i(f))}"] \arrow[d, "\cong"] & \H_n(H; \K_i(B)) \arrow[d, "\cong"] \\
 \bb L_n \K_i(G \ltimes A) \arrow[r, "{\bb L_n (\Ind_\Omega, \alpha_\Omega, f)_i}"] & \bb L_n \K_i(H \ltimes B)
\end{tikzcd} \]
\begin{proof}
For $n \geq 0$ let $P_n = (\Ind^G_{G^0} \Res^{G^0}_G)^{n+1} A$ and $Q_n = (\Ind^H_{H^0} \Res^{H^0}_H)^{n+1} B$ and consider the projective resolutions $P_\bullet \to A$ and $Q_\bullet \to B$ induced by the adjunctions as in \Cref{adjunctions are good}. Recall that $\bb L_n (\Ind_\Omega, \alpha_\Omega, f)$ arises from a chain map $\tilde f \colon P_\bullet \to \Ind_\Omega Q_\bullet$ over $f \colon A \to \Ind_\Omega B$ as the map induced in homology by the composite chain map $\alpha_\Omega(Q) \circ \K_i(G \ltimes \tilde f)$:
\[\begin{tikzcd}
	\cdots & {\K_i(G \ltimes P_n )} & \cdots & {\K_i(G \ltimes P_0 )} & {} \\
	\cdots & {\K_i(G \ltimes \Ind_\Omega Q_n)} & \cdots & {\K_i(G \ltimes \Ind_\Omega Q_0 )} \\
	\cdots & {\K_i( H \ltimes Q_n )} & \cdots & {\K_i(H \ltimes Q_0 )}
	\arrow[from=1-2, to=2-2, "\K_i(G \ltimes \tilde f_n)"]
	\arrow[from=2-2, to=3-2, "\alpha_\Omega(Q_n)"]
	\arrow[from=1-2, to=1-3]
	\arrow[from=1-3, to=1-4]
	\arrow[from=2-2, to=2-3]
	\arrow[from=2-3, to=2-4]
	\arrow[from=1-4, to=2-4, "\K_i(G \ltimes \tilde f_0)"]
	\arrow[from=2-4, to=3-4, "\alpha_\Omega(Q_0)"]
	\arrow[from=3-2, to=3-3]
	\arrow[from=3-3, to=3-4]
	\arrow[from=2-1, to=2-2]
	\arrow[from=1-1, to=1-2]
	\arrow[from=3-1, to=3-2]
\end{tikzcd}\]
On the module side, the modules $\K_i(P_n)$ are left $G$-acyclic and $\K_i(Q_n)$ are left $H$-acylic. Furthermore, the chain complexes $\K_i(P_\bullet)$ and $\K_i(\Ind_\Omega Q_\bullet)$ are exact because $\ker \K_i \supseteq \I^G_0$ and similarly for $H$, and so form acyclic resolutions $\K_i(P_\bullet) \to \K_i(A)$ and $\K_i(Q_\bullet) \to \K_i(B)$. The composite chain map $\zeta_{\Omega,Q}^{-1} \circ \K_i(\tilde f)$ shown below lifts $\K_i(f) \colon \K_i(A) \to \Ind_\Omega \K_i(B)$, noting that $\K_i(f)$ is really shorthand for $\zeta_{\Omega,B}^{-1} \circ \K_i(f)$.
\[\begin{tikzcd}
	\cdots & {\K_i(P_n )} & \cdots & {\K_i(P_0)} & {} \\
	\cdots & {\K_i(\Ind_\Omega Q_n )} & \cdots & {\K_i(\Ind_\Omega Q_0 )} \\
	\cdots & {\Ind_\Omega \K_i(Q_n)} & \cdots & {\Ind_\Omega \K_i(Q_0)}
	\arrow[from=1-2, to=2-2, "\K_i(\tilde f_n)"]
	\arrow[from=1-2, to=1-3]
	\arrow[from=1-3, to=1-4]
	\arrow[from=2-2, to=2-3]
	\arrow[from=2-3, to=2-4]
	\arrow[from=1-4, to=2-4, "\K_i(\tilde f_0)"]
	\arrow[from=2-1, to=2-2]
	\arrow[from=1-1, to=1-2]
	\arrow[from=3-1, to=3-2]
	\arrow[from=3-2, to=3-3]
	\arrow[from=3-3, to=3-4]
	\arrow[from=2-2, to=3-2, "\zeta_{\Omega, Q_n}^{-1}"]
	\arrow[from=2-4, to=3-4, "\zeta_{\Omega,Q_0}^{-1}"]
\end{tikzcd}\] 
Therefore by \cite[Theorem 3.5]{Miller24}, the map $\H_*(\Omega;\K_i(f))$ is induced by the chain map $(\delta_\Omega \otimes \id) \circ (\zeta_{\Omega,Q_n}^{-1} \circ \K_i(\tilde f_n))_G$. By \Cref{modules and crossed product comparison} the diagram
\[ 
\begin{tikzcd}[column sep = 4cm]
\K_i(P_n)_G \arrow[d, "\cong", "\gamma_{P_n}"'] \arrow[r, "(\delta_\Omega \otimes \id) \circ (\zeta_{\Omega,Q_n}^{-1} \circ \K_i(\tilde f_n))_G"] & \K_i(Q_n)_H \arrow[d, "\cong", "\gamma_{Q_n}"'] \\
\K_i(G \ltimes P_n) \arrow[r, "\alpha_\Omega(Q_n) \circ \K_i(G \ltimes \tilde f_n )"] & \K_i(H \ltimes Q_n) 
\end{tikzcd}
\]
commutes for each $n \geq 0$. On the level of homology the vertical chain maps induce the identification in \Cref{homology identified as a derived functor}, and the horizontal chain maps induce $\H_*(\Omega; \K_i(f))$ and $\bb L_*(\Ind_\Omega, \alpha_\Omega, f)_i$ respectively, so we are done.
\end{proof}
\end{thm}

\section{Functoriality of the ABC spectral sequence}\label{ABC functor section}

In this section we prove the following functorial version of Theorem \ref{ABC}.

\begin{thm}\label{ABCthmintext}
An ABC morphism $\fk m \colon \fk M \to \fk M'$ induces functorially a morphism of ABC spectral sequences $\ABC(\fk m) \colon \ABC(\fk M) \to \ABC(\fk M')$, such that:
\begin{enumerate}[label=(\roman*)]
\item the map on the second sheet is given by $\bb L_p(\fk m)_q \colon \bb L_p(\fk M)_q \to \bb L_p (\fk M')_q$, 
\item the map on the limit sheet agrees with $\bb L(\fk m) \colon \bb L (\fk M) \to \bb L (\fk M')$.
\end{enumerate}
\end{thm}

We first review Meyer's construction of the ABC spectral sequence associated to an ABC tuple $(\T,\I,F,A)$ from \cite[Section 4]{Meyer08}. This starts by picking an $\I$-projective resolution over $A$ and embedding this in a phantom tower. Note the convention that a degree one map $f \colon A \to \Sigma B$ is written $f \colon A \circledrightarrow B$.
\begin{defn}[Phantom tower]
Let $A$ be an object in a triangulated category $\T$ with a homological ideal $\I$. A \textit{pre-$\I$-phantom tower} over $A$ is a diagram
\[\begin{tikzcd}[column sep=small]
A \arrow[r,  "="{description}, draw=none] &  N_0 \arrow[rr, "\iota_0^1"] &                         & N_1 \arrow[ld, circled, "\epsilon_0"] \arrow[rr, "\iota_1^2"] &                                                        & N_2 \arrow[ld, circled, "\epsilon_1"] \arrow[rr] &                                                        & \cdots \arrow[ld, circled] &                                      \\
                                 &                             & P_0 \arrow[lu, "\pi_0"] &                                                              & P_1 \arrow[ll, circled, "\delta_1"] \arrow[lu, "\pi_1"] &                                                 & P_2 \arrow[ll, circled, "\delta_2"] \arrow[lu, "\pi_2"] &                           & \cdots \arrow[ll, circled] \arrow[lu]
\end{tikzcd}\]
written $\mathcal P \to A$ with $\mathcal P = (P_\bullet, N_\bullet,\delta,\pi,\iota,\epsilon)$ or $(P_\bullet, N_\bullet)$ in brief, such that:
\begin{itemize}
\item The top triangles are exact, i.e. $P_n \xrightarrow{\pi_n} N_n \xrightarrow{\iota^{n+1}_n} N_{n+1} \xcircledrightarrow{\epsilon_n} P_n $ is an exact triangle for each $n$.
\item The bottom triangles commute, i.e. $\delta_n = \epsilon_{n-1} \circ \pi_n$ for each $n$.
\item The top maps are $\I$-phantom, i.e. $\iota^{n+1}_n \in \I$ for each $n$.
\item The sequence  $\cdots \xcircledrightarrow{\delta_2} P_1 \xcircledrightarrow{\delta_1} P_0 \xrightarrow{\pi_0} N_0 \rightarrow 0 $ is $\I$-exact. 
\end{itemize}
The last two conditions are equivalent in the presence of the other conditions. A \textit{tower map} $\mathcal P \to \mathcal Q$ of pre-$\I$-phantom towers $\mathcal P \to A$ and $\mathcal Q \to B$ over a morphism $f \colon A \to B$ is a morphism of diagrams compatible with $f$. We say that $\mathcal P$ is \textit{$\I$-phantom} if $P_n$ is $\I$-projective for each $n$. 
\end{defn}

Meyer shows that any $\I$-projective resolution embeds in an $\I$-phantom tower (note that the suspension of an $\I$-projective object is $\I$-projective) and that any morphism lifts to a map of phantom towers over the objects \cite[Lemma 3.3]{Meyer08}. We must additionally consider pre-phantom towers, but the same proof\footnote{using the slightly modified version \Cref{resolutionmap} of \cite[Proposition 3.26]{MeyNes10}} gives the following statement:

\begin{lem}\label{towermap}
Any $\I$-projective resolution $P_\bullet \rightarrow A$ can be embedded in an $\I$-phantom tower, which is unique up to non-canonical isomorphism.

Furthermore, if $ \cal P $ is a pre-$\I$-phantom tower over $A$, $\cal Q$ is a pre-$\I$-phantom tower over $B$, and $P_\bullet$ and $Q_\bullet$ are the $\I$-resolutions of $A$ and $B$ embedded in $\cal P$ and $\cal Q$ respectively, then any chain map from $P_\bullet$ to $Q_\bullet$ over a map $f\colon A \to B$ extends to a tower map $\cal P \rightarrow \cal Q$ over $f$. In particular, if $\cal P$ is $\I$-phantom, any map $f\colon A \rightarrow B$ lifts (non-canonically) to a tower map $\cal P \rightarrow \cal Q$.
\end{lem}

Given a pre-$\I$-phantom tower $\mathcal P \to A$ we construct an exact couple $\EC(\mathcal P)$ as follows. The top collection of exact triangles in $\mathcal P$ forms a single triangle of $\Z$-graded objects after defining $N_n = A$, $\iota^{n+1}_n = \id_A$ and $P_n = 0$ for $n < 0$.
\[ \begin{tikzcd}[column sep=small]
N_\bullet \arrow[rr, "\iota"] &                     & N_\bullet \arrow[ld, circled, "\epsilon"] \\
                      & P_\bullet \arrow[lu, "\pi"] &                                 
\end{tikzcd} \]
We now consider the stable homological functor $F \colon  \T \to \cat{Ab}_*$. We use this to construct an exact triangle of bigraded abelian groups
\begin{equation}\label{exact couple construction}
\begin{tikzcd}[column sep=small]
D \arrow[rr, "i"] &                     & D \arrow[ld, "j"] & D_{p,q} = F_{p+q+1}(N_{p+1}) \\
                      & C \arrow[lu, "k"] &                                 & C_{p,q} = F_{p+q}(P_p)
\end{tikzcd} 
\end{equation}
with maps $i$, $j$ and $k$ defined at $p,q \in \Z$ by\footnote{The notation $F_n(f)$ for the map $F_n(A) \to F_{n-1}(B)$ induced by a degree one morphism $f \colon A \circledrightarrow B$ may not be standard.} 
\begin{align*}
 i_{p,q} & \defeq F_{p+q+1}(\iota_{p+1}^{p+2})  \colon  D_{p,q}  \rightarrow D_{p+1,q-1} & \text{deg } i & = (1,-1) \\
 j_{p,q} & \defeq F_{p+q+1}(\epsilon_p)  \colon  D_{p,q}  \rightarrow C_{p,q} & \text{deg } j & = (0,0) \\
 k_{p,q} & \defeq F_{p+q}(\pi_p)  \colon  C_{p,q}  \rightarrow D_{p-1,q} & \text{deg } k & = (-1,0).
\end{align*}
This forms an exact couple $\EC(\mathcal P) = (C,D,i,j,k)$ which is functorial with respect to tower maps. The construction of a spectral sequence from an exact couple is standard, see for example \cite[p336--337]{Maclane95}, and yields a spectral sequence $(E,d)$ starting from the first sheet. We list some important properties:
\begin{itemize}
\item The first sheet $E^1$ satisfies $E^1_{p,q} = F_{p+q}(P_p)$ for $p \geq 0$, $q \in \Z$ and vanishes for $p \leq -1$. This implies $E^r_{p,q} = 0$ for any $r \geq 1$ and $q \in \Z$ when $p \leq -1$.
\item The first differential $d^{(1)}_{p,q} \colon E^1_{p,q}  \to E^1_{p-1,q} $ is given by $F_{p+q}(\delta_p) \colon F_{p+q}(P_p) \to F_{p+q-1}(P_{p-1})$ for $p \geq 0$, $q \in \Z$.
\item The second sheet $E^2 = \frac{k^{-1} ( \im i )}{j( \ker i)}$ is given at $p \geq 0$, $q \in \Z$ by $\frac{\ker F_{p+q}(\delta_p)}{\im F_{p+q+1}(\delta_{p+1})}$.
\item For $r \geq 1$ the $r$th sheet $E^r = \frac{k^{-1}(\im i^{r-1})}{j(\ker i^{r-1})}$ has differential $d^{(r)} \colon E^r \to E^r$ satisfying $d^{(r)}(x) = j(y)$ for $x \in E^r$ and $y \in D$ with $i^{r-1}(y) = k(x)$.
\item The limit sheet $E^\infty$ is $\frac{\bigcap_{r \geq 1} k^{-1} ( \im i^{r} )}{ \bigcup_{r \geq 1} j ( \ker i^{r} )}$.
\end{itemize}
The \textit{ABC spectral sequence} $\ABC(\T,\I,F,A)$ is the spectral sequence $(E^r,d^{(r)})_{r \geq 2}$ starting from the second sheet, for any $\I$-phantom tower $\mathcal P \to A$. In this case, the description of the second sheet gives an isomorphism $E^2_{p,q} \cong \bb L_{p} F_q(A)$. For any other $\I$-phantom tower $\mathcal Q \to A$, there is a tower map $\mathcal P \to \mathcal Q$ over the identity $\id \colon A \to A$ by \Cref{towermap}. This induces a morphism of exact couples $\EC(\mathcal P) \to \EC(\mathcal Q)$ and therefore a morphism of spectral sequences. On the second sheets, this is induced by a chain map of $\I$-projective resolutions $P_\bullet \to Q_\bullet$ of $A$, and thus is the identity on $E^2_{p,q} = \bb L_p F_q(A)$ by \Cref{resolutionmap}. Since the later sheets are subquotients of the second sheet, the spectral sequences induced by $\mathcal P$ and $\mathcal Q$ are identical from the second sheet onwards.

\begin{lem}
Let $(\Phi, \alpha, f) \colon (\T,\I, F, A) \to (\T', \I', F', A')$ be an ABC cycle and let $\cal P'$ be an $\I'$-phantom tower over $A'$. Then $\Phi(\cal P')$ is a pre-$\I$-phantom tower over $\Phi(A')$.
\begin{proof}
The triangulated functor $\Phi$ preserves exact triangles and sends $\I'$ into $\I$, and thus also sends $\I'$-exact sequences to $\I$-exact sequences.
\end{proof}
\end{lem}

\begin{prop}\label{morphism of ABC spectral sequences}
An ABC morphism $\fk m = [\Phi, \alpha, f] \colon (\T,\I, F, A) \to (\T', \I', F', A')$ functorially induces a morphism 
\[\ABC(\fk m) \colon \ABC(\T, \I, F, A) \to \ABC(\T', \I', F', A')\]
of spectral sequences with the following property. Let $\cal P$ be an $\I$-phantom tower over $A$, let $\cal P'$ be an $\I'$-phantom tower over $A'$, and let $\tilde{f} \colon \cal P \to \Phi(\cal P')$ be a tower map over $f \colon A \to \Phi(A')$. Let $\alpha_{\cal P'} \colon \EC(\Phi(\cal P')) \to \EC(\cal P')$ be the morphism of exact couples induced by the natural transformation $\alpha$. Then the morphism of spectral sequences induced by $\alpha_{\cal P'} \circ \EC(\tilde{f}) \colon \EC(\cal P) \to \EC(\cal P')$ agrees with $\ABC(\fk m)$ from the second sheet onwards. On the second sheets, we have $\ABC(\fk m)^2_{p,q} = \bb L_p(\fk m)_q \colon \bb L_p F_q(A) \to \bb L_p F'_q(A')$.
\[ \begin{tikzcd}[column sep = 0.2, row sep = 0.2]
 N_0 \arrow[rr] \arrow[ddd, red] &                & N_1 \arrow[ddd, red] \arrow[ldd, circled] \arrow[rr] &                                   & \phantom{l} \cdots \phantom{l} \arrow[ldd, circled] &                                      & F_*(N_\bullet) \arrow[rr] &    & F_*(N_\bullet) \arrow[ldd]      \\
  & & & & & & & & {\phantom{A}} \\
             {}             & P_0  \arrow[luu] & {}                                & P_1  \arrow[ll, circled, crossing over] \arrow[luu] & {}                        & \phantom{l} \cdots \phantom{l} \arrow[ll, circled] \arrow[luu] &  & F_{*}(P_\bullet) \arrow[luu]  &      \\
 \Phi(N_0') \arrow[rr] &                & \Phi(N_1') \arrow[ldd, circled] \arrow[rr] &                                   & \phantom{l} \cdots \phantom{l} \arrow[ldd, circled] &                                      & F_* \Phi(N'_\bullet) \arrow[rr] \arrow[from = uuu, red] &    & F_* \Phi(N'_\bullet) \arrow[ldd] \arrow[from = uuu, red]   \\
   & & & & & & & & {\phantom{A}} \\
                          & \Phi(P_0') \arrow[from=uuu, red, crossing over] \arrow[luu] &                                   & \Phi(P_1') \arrow[from=uuu, red, crossing over] \arrow[ll, circled] \arrow[luu] &                           & \phantom{l} \cdots \phantom{l}\arrow[ll, circled] \arrow[luu] &  & {F_{*} \Phi(P'_\bullet)} \arrow[luu] \arrow[from = uuu, red, crossing over] &  \\
 & & & & &                                      & F'_* (N'_\bullet) \arrow[rr] \arrow[from = uuu, red] &    & F'_* (N'_\bullet) \arrow[ldd] \arrow[from = uuu, red]     \\
   & & & & & & & &   {\phantom{A}}    \\
        & & & & & &  & {F'_{*} (P'_\bullet)} \arrow[luu] \arrow[from = uuu, red, crossing over] &   
\end{tikzcd}  \]
In red on the left is the tower map $\tilde{f} \colon \cal P \to \Phi ( \cal P')$, and on the right are the morphisms $\EC(\tilde{f})$ and $\alpha_{\cal P'}$ of exact couples.
\begin{proof}
To construct $\ABC(\fk m)$, we consider arbitrary phantom towers $\cal P$ and $\cal P'$ over $A$ and $A'$ respectively. By \Cref{towermap}, there is a tower map $\tilde{f} \colon \cal P \to \Phi(\cal P')$ over $f \colon A \to \Phi(A')$. We obtain a morphism of exact couples $ \alpha_{\cal P'} \circ EC(\tilde{f})$. We set $\ABC(\fk m)$ to be the morphism of spectral sequences induced by this morphism, starting at the second sheet.

Let $\hat{f} \colon P_\bullet \to \Phi(P'_\bullet)$ be the restriction of $\tilde{f}$ to the $\I$-exact sequence $P_\bullet$ inside $\cal P$. By the description of the exact couple to spectral sequence functor, the map $\ABC(\fk m)^2_{p,q} \colon \bb L_p F_q(A) \to \bb L_p F_q(A')$ is induced in homology by the chain map $\alpha_{P'} \circ F(\hat{f})$. \Cref{derivedmaps} then tells us that $\ABC(\fk m)^2_{p,q} = \bb L_p (\fk m)_q$. Morphisms of spectral sequences are determined by their restrictions to the second sheets (see \cite[Mapping Lemma 5.2.4]{Weibel94}), and therefore $\ABC(\fk m)$ is functorial and independent of the choice of $\I$-phantom tower $\cal P$ and tower map $\tilde{f} \colon \cal P \to \Phi(\cal P')$.
\end{proof}
\end{prop}

Constructing functoriality of spectral sequences alone is not the full story; we care about the convergence of spectral sequences, which is extra structure. Recall that a spectral sequence $E$ converges (strongly) to a graded abelian group $G$ with (ascending) filtration $(\cal F_k G)_{k \geq 0}$ with $\mathcal F_0 G = 0$ if the filtration is exhaustive ($\bigcup_{k \geq 0} \mathcal F_k G = G$) and for each $p,q \in \Z$ there is an isomorphism
\[E^\infty_{p,q} \cong \frac{\cal F_{p+1} G_{p+q}}{\cal F_p G_{p+q}}. \]
The ABC spectral sequence of an ABC tuple $(\T,\I,F,A)$ converges strongly to the graded abelian group $\bb L F(A)$ with the \textit{phantom filtration} $(\bb L F \colon \I^k (A))_{k \geq 0}$ \cite[Theorem 5.1]{Meyer08}:
\begin{defn}[Phantom filtration]
Let $\I$ be a homological ideal in a triangulated category $\T$, let $F\colon  \T \to \Ab_*$ be a stable homological functor and let $A$ be an object of $\T$. For each $k \geq 0$, set
\[ F \colon  \I^k (A) \defeq \{ x \in F(A) \mid  F(f)(x) = 0 \; \text{for all} \; f \in \I^k(A,B), \; B \in \T \}. \]
This forms an ascending filtration of $F(A)$ called the \textit{$\I$-phantom filtration} of $F$ at $A$.
\[ 0  = F \colon  \I^0 (A) \subseteq F \colon  \I^1 (A) \subseteq F \colon  \I^2 (A) \subseteq \cdots \subseteq F(A)  \]
\end{defn}
To show the convergence of the ABC spectral sequence of the ABC tuple $\fk M = (\T,\I,F,A)$, Meyer shows that $(\bb L F \colon \I^k (A))_{k \geq 0}$ is exhaustive and constructs isomorphisms
\[ \psi_{p,q}(\fk M) \colon \ABC(\fk M)^\infty_{p,q} \to \frac{\bb L F_{p+q} \colon \I^{p+1}(A)}{\bb L F_{p+q} \colon \I^{p}(A)}. \]
Our aim is to show that these isomorphisms are functorial with respect to ABC morphisms. We collect the following description of the phantom filtration:
\begin{lem}\label{kerneldescription}
Let $(\T,\I,F,A)$ be an ABC tuple with associated localisation functor $L \colon \T \to \T$ and let $\mathcal P \to A$ be an $\I$-phantom tower over $A$:
\[\begin{tikzcd}[column sep=small]
A \arrow[r,  "="{description}, draw=none] & N_0 \arrow[rr, "\iota_0^1"] &                         & N_1 \arrow[ld, circled, "\epsilon_0"] \arrow[rr, "\iota_1^2"] &                                                        & N_2 \arrow[ld, circled, "\epsilon_1"] \arrow[rr] &                                                        & \cdots \arrow[ld, circled] &                                      \\
                                 &                             & P_0 \arrow[lu, "\pi_0"] &                                                              & P_1 \arrow[ll, circled, "\delta_1"] \arrow[lu, "\pi_1"] &                                                 & P_2 \arrow[ll, circled, "\delta_2"] \arrow[lu, "\pi_2"] &                           & \cdots \arrow[ll, circled] \arrow[lu]
\end{tikzcd}\]
Then $\ker F(\iota^k_0) = F \colon  \I^{k} (A)$ for every $k \geq 0$, where $\iota^k_0 \defeq \iota^k_{k-1} \circ \iota^{k-1}_{k-2} \circ \cdots \circ \iota^{1}_0$. Moreover, the subgroups $\bb L F \colon  \I^k (A)$ and $ F \colon  \I^k (LA)$ of $\bb L F (A)$ are equal.
\begin{proof}
That $\ker F(\iota^k_0) = F \colon  \I^{k} (A)$ follows immediately from \cite[Lemmas 3.5-3.7]{Meyer08}. The localisation functor $L$ maps $\I$-phantom towers to $\I$-phantom towers, so we may apply the first part of the lemma to $(\bb L F, \mathcal P)$ and $(F, L(\mathcal P))$ respectively to conclude $\bb L F \colon  \I^k (A) = \ker F(L(\iota_0^k)) =
F \colon \I^k(LA)$.
\end{proof}
\end{lem}

ABC morphisms respect the phantom filtrations in the following sense:
\begin{lem}\label{application respects filtration}
Let $[\Phi, \alpha, f] \colon (\T,\I,F,A) \to (\T', \I', F', A')$ be an ABC morphism. Then $\App(\Phi, \alpha, f) (F \colon \I^k (A)) \subseteq  F' \colon {\I'}^k (A')$.
\begin{proof}
Let $x \in F \colon \I^k (A)$, and suppose $g \in {\I'}^k(A',B')$ for some $B' \in \T'$. To show that $\App(\Phi, \alpha, f)(x) \in F' \colon {\I'}^k (A')$ we must show that $F'(g) \circ \alpha_{A'} \circ F(f) (x)$ vanishes. By naturality of $\alpha$ at $g$ we may rewrite this as $\alpha_{B'} \circ F(\Phi(g) \circ f)(x)$, which vanishes because $\Phi(g) \circ f \in \I^k$.
\end{proof}
\end{lem}
In particular, it follows that $\bb L (\Phi, \alpha, f) \colon \bb L F (A) \to \bb L F'(A')$ respects the phantom filtrations. The following lemma will help us to describe the convergence of the ABC spectral sequence:
\begin{lem}\label{map from infinity sheet}
Let $(C,D,i,j,k)$ be an exact couple, let $(E,d)$ be the associated spectral sequence and let $p,q \in \Z$ with $p \geq 0$. Recall that we have the following descriptions.
\[ \begin{tikzcd}[column sep = small]
D \arrow[rr, "{(1,-1)}", "i"'] &                                   & D \arrow[ld, "{(0,0)}", "j"'] \\
                                  & C \arrow[lu, "{(-1,0)}", "k"'] &                                 
\end{tikzcd} \qquad E^{\infty}_{p,q} = \frac{\bigcap_{r \geq 1}(k^{-1}(\im i^r))_{p,q}}{\bigcup_{r \geq 1}(j ( \ker i^r ))_{p,q}}.  \]
Then there is a homomorphism $\psi_{p,q} \colon E^\infty_{p,q} \to \frac{(\ker i^{p+1})_{-1,p+q}}{(\ker i^{p})_{-1,p+q}} $ which makes the following diagram commute:
\[\begin{tikzcd}
	{\bigcap_{r \geq 1} k^{-1}(\im i^r)_{p,q}} & {D_{p-1,q}} \\
	{E^\infty_{p,q}} & {\frac{(\ker i^{p+1})_{-1,p+q}}{(\ker i^{p})_{-1,p+q}}}
	\arrow["k", from=1-1, to=1-2]
	\arrow["\psi_{p,q}", from=2-1, to=2-2]
	\arrow[two heads, from=1-1, to=2-1]
	\arrow["{i^p}", hook, from=2-2, to=1-2]
\end{tikzcd}\]
\begin{proof}
The map $k \colon \bigcap_{r \geq 1} k^{-1}(\im i^r)_{p,q} \to D_{p-1,q}$ factors through $E^\infty_{p,q}$ because $k \circ j = 0$. We need only check that the image is contained in $i^p((\ker i^{p+1})_{-1,p+q})$. For $x \in \bigcap_{r \geq 1} k^{-1}(\im i^r)_{p,q}$ there is $y \in D_{-1,p+q}$ such that $k(x) = i^p(y)$ by construction. Because $i \circ k = 0$, it follows that $y \in \ker i^{p+1}$.
\end{proof}
\end{lem}
For an ABC tuple $\fk M = (\T,\I,F,A)$, \Cref{kerneldescription} tells us that for the exact couple $\EC(\mathcal P)$ associated to any $\I$-phantom tower $\mathcal P \to A$ and $p,q \in \Z$ with $p \geq 0$ we have
\[\frac{(\ker i^{p+1})_{-1,p+q}}{(\ker i^{p})_{-1,p+q}} = \frac{F_{p+q} \colon \I^{p+1} (A)}{F_{p+q} \colon \I^{p} (A)}.\]
\Cref{map from infinity sheet} gives us the map 
\[ \varphi_{p,q}(\fk M) \colon \ABC(\fk M)^\infty_{p,q} \to \frac{F_{p+q} \colon \I^{p+1} (A)}{F_{p+q} \colon \I^{p} (A)}\]
appearing in \cite[Proposition 4.3]{Meyer08}. Meyer notes that this does not depend on the choice of $\I$-phantom tower and shows that it is an isomorphism whenever $A \in \fk P_\I $ \cite[Proposition 4.10]{Meyer08}. This applies in particular to $\fk L(\fk M) \defeq (\T,\I,F,LA)$, thus $\varphi_{p,q}(\fk L(\fk M))$ is an isomorphism. The $\I$-equivalence $\mu_A \colon LA \to A$ induces an isomorphism $\ABC(\mu_A)$ of ABC spectral sequences because the restriction $\bb L F (\mu_A)$ to the second sheets is an isomorphism. Putting these together forms the desired isomorphism $\psi(\fk M)_{p,q}$:
\begin{equation}\label{convergence map psi}
\begin{tikzcd}[column sep = 5em]
	{\ABC(\fk M)^\infty_{p,q}} & {\frac{\bb L F_{p+q} \colon \I^{p+1}(A)}{\bb L F_{p+q} \colon \I^{p}(A)}} \\
	{\ABC(\fk L(\fk M))^\infty_{p,q}} & {\frac{F_{p+q} \colon \I^{p+1}(L A)}{F_{p+q} \colon \I^{p}(L A)}}
	\arrow["{\psi(\fk M)_{p,q}}", from=1-1, to=1-2]
	\arrow["{\varphi(\fk L(\fk M))_{p,q}}", from=2-1, to=2-2]
	\arrow["{\ABC(\mu_A)^\infty_{p,q}}", "\cong"', from=2-1, to=1-1]
	\arrow["\rotatebox{90}{$\,=$}"{description}, phantom, from=2-2, to=1-2]
\end{tikzcd}
\end{equation}
\begin{prop}\label{naturality of limit sheet phi prop}
Let $\fk M = (\T,\I,F,A)$ and $\fk M' = (\T',\I',F',A')$ be ABC tuples and let $\fk m \colon \fk M \to \fk M'$ be an ABC morphism. Then the following diagram commutes:
\begin{equation}\label{naturality of limit sheet phi}
 \begin{tikzcd}[column sep = 5em]
\ABC(\fk M)^{\infty}_{p,q} \arrow[r, "\varphi(\fk M)_{p,q}"] \arrow[d, "\ABC(\fk m)^{\infty}_{p,q}"]  & \frac{F_{p+q} \colon  \I ^{p+1} (A)}{F_{p+q} \colon  \I ^{p} (A)} \arrow[d, " \App(\fk m)_{p,q} "] \\
{\ABC(\fk M')}^{\infty}_{p,q} \arrow[r, "\varphi(\fk M')_{p,q}"] & \frac{F'_{p+q} \colon  {\I'} ^{p+1} (A')}{F'_{p+q} \colon  {\I'} ^{p} (A')}
\end{tikzcd} 
\end{equation}
\begin{proof}
Let $\mathcal P \to A$ and $\mathcal P' \to A'$ be phantom towers over $A$ and $A'$ and let $\EC(\mathcal P) = (C,D,i,j,k)$ and $\EC(\mathcal Q) = (C',D',i',j',k')$ be the associated exact couples. Suppose $\fk m = [\Phi, \alpha, f]$ and let $\tilde{f} \colon \cal P \to \Phi(\cal P')$ be a tower map over $f \colon A \to \Phi(A')$. Let $S \colon C \to C'$ and $T \colon D \to D'$ be the constituent components of the morphism $\alpha_{\cal P'} \circ \EC(\tilde{f})$ of exact couples. Then for each $n \in \Z$ we have $D_{-1,n} = F_n(A)$ and $D'_{-1,n} = F'_n(A')$. Furthermore, $T_{-1,n} \colon D_{-1,n} \to D'_{-1,n}$ is given by the map $\App(\fk m)_n \colon F_n(A) \to F'_n(A')$.
\[\begin{tikzcd}[column sep = small, row sep = small]
D \arrow[dd, "T", red] \arrow[rr, "i"] &                                   & D \arrow[dd, "T", red] \arrow[ld, "j"] &  & {D_{-1,n}} \arrow[r, "="{description}, phantom] \arrow[dd, "{T_{-1,n}}", red] & F_n(A) \arrow[dd, "\App(\fk m)_n", red] \\
                                  & C \arrow[lu, "k"]  &                                   &  &                                                    &                                    \\
D' \arrow[rr, "i'", pos = 0.3]               &                                   & D' \arrow[ld, "j'"]               &  & {D'_{-1,n}} \arrow[r, "="{description}, phantom]                         & F'_n(A')                           \\
                                  & C' \arrow[lu, "k'"]   \arrow[from = uu, crossing over, "S", pos = 0.3, red]            &                                   &  &                                                    &                                   
\end{tikzcd}\]
Recall that $\ABC(\fk M)^\infty_{p,q} = \frac{\bigcap_{r \geq 1}(k^{-1}(\im i^r))_{p,q}}{\bigcup_{r \geq 1}(j ( \ker i^r ))_{p,q}}$, let $x \in \bigcap_{r \geq 1}(k^{-1}(\im i^r))_{p,q}$ and pick $y \in F_{p+q} \colon \I^{p+1}(A)$ with $\varphi(\fk M)_{p,q} ([x]) = [y]$. Since $y \in D_{-1,p+q} = F_{p+q}(A)$, we may use that $\App(\fk m)_{p+q} = T_{-1,p+q}$ to calculate that 
\[\App(\fk m)_{p,q}(\varphi(\fk M)_{p,q}([x])) = \App(\fk m)_{p,q}([y]) = \left[ \App(\fk m)_{p+q}(y) \right] = [T(y)].\]
On the other hand,  $\ABC(\fk m)^\infty_{p,q} ([x]) = [S(x)]$. We may calculate that ${i'}^p (T(y)) = T( i^p(y)) = T( k(x)) = k'(S(x))$, and therefore by \Cref{map from infinity sheet}, $\varphi(\fk M')_{p,q}([S(x)]) = [T(y)]$. We conclude that \eqref{naturality of limit sheet phi} commutes.
\end{proof}
\end{prop}

\begin{thm*}[Theorem \ref{ABCthmintext}]
An ABC morphism $\fk m \colon \fk M \to \fk M'$ functorially induces a morphism of ABC spectral sequences $\ABC(\fk m) \colon \ABC(\fk M) \to \ABC(\fk M')$, such that:
\begin{enumerate}[label=(\roman*)]
\item the map on the second sheets is given by $\bb L_p(\fk m)_q \colon \bb L_p(\fk M)_q \to \bb L_p (\fk M')_q$, 
\item the map on the limit sheets agrees with $\bb L(\fk m) \colon \bb L (\fk M) \to \bb L (\fk M')$.
\end{enumerate}
\begin{proof}
Let $\fk m \colon \fk M \to \fk M'$ be represented by an ABC cycle
\[(\Phi, \alpha, f) \colon (\T, \I, F, A) \to (\T', \I', F', A').\]
The construction of the morphism $\ABC(\fk m) \colon \ABC(\fk M) \to \ABC(\fk M')$ was given in \Cref{morphism of ABC spectral sequences}, where it was shown to be functorial and agree with $\bb L_p(\fk m)_q$ on the second sheets. To show compatibility of the map on the limit sheets with $\bb L(\fk m)$ we must show that
\[ \begin{tikzcd}[column sep = 5em]
	{\ABC(\fk M)^{\infty}_{p,q}} & {\frac{\bb L F_{p+q} \colon \I^{p+1}(A)}{\bb L F_{p+q} \colon \I^{p}(A)}} \\
	{\ABC(\fk M')^{\infty}_{p,q}} & {\frac{\bb L' F'_{p+q} \colon {\I'}^{p+1}(A')}{\bb L' F'_{p+q} \colon {\I'}^{p}(A')}}
	\arrow[from=1-1, to=1-2, "\cong", "\psi(\fk M)_{p,q}"']
	\arrow[from=2-1, to=2-2, "\cong", "\psi(\fk M')_{p,q}"']
	\arrow[from=1-1, to=2-1, "\ABC(\fk m)^\infty_{p,q}"]
	\arrow[from=1-2, to=2-2, "\bb L(\fk m)_{p,q}"]
\end{tikzcd} \]
commutes. Let $\mu \colon L \Rightarrow \id$ and $\mu' \colon L' \Rightarrow \id$ be the localisation natural transformations associated to $\fk M$ and $\fk M'$ and set $\fk L(\fk M) = (\T,\I,F,LA)$, $\fk L(\fk M') = (\T', \I', F', L'A')$ and $\fk L(\fk m) = [\Phi, \alpha, \bb L_\Phi(A') \circ L(f)]$. By \eqref{convergence map psi} the diagram breaks into three parts:
\begin{equation*}\label{breakdown CD}
\begin{tikzcd}[column sep = small]
	{\ABC(\fk M)^{\infty}_{p,q}} & & {\ABC(\fk L(\fk M))^{\infty}_{p,q}} & & {\frac{F_{p+q} \colon \I^{p+1}(L A)}{F_{p+q} \colon \I^{p}(L A)}} & {\frac{\bb L F_{p+q} \colon \I^{p+1}(A)}{\bb L F_{p+q} \colon \I^{p}(A)}} \\
	{\ABC(\fk M')^{\infty}_{p,q}} & & {\ABC(\fk L(\fk M'))^{\infty}_{p,q}} & & {\frac{F'_{p+q} \colon {\I'}^{p+1}(L' A')}{F'_{p+q} \colon {\I'}^{p}(L' A')}} & {\frac{\bb L' F'_{p+q} \colon {\I'}^{p+1}(A')}{\bb L' F'_{p+q} \colon {\I'}^{p}(A')}}
	\arrow[from=1-1, to=1-3, "\cong"]
	\arrow[from=1-3, to=1-5, "\cong", "\varphi(\fk M)_{p,q}"']
	\arrow[from=1-5, to=1-6, "="{description}, phantom]
	\arrow[from=1-6, to=2-6, "\bb L ( \fk m)_{p,q}"]
	\arrow[from=2-5, to=2-6, "="{description}, phantom]
	\arrow[from=1-5, to=2-5, "\App(\fk L(\fk m))_{p,q}"]
	\arrow[from=2-3, to=2-5, "\cong", "\varphi(\fk M)_{p,q}"']
	\arrow[from=1-3, to=2-3, "\ABC(\fk L (\fk m))^\infty_{p,q}"]
	\arrow[from=1-1, to=2-1, "\ABC( \fk m )^{\infty} _{p,q} "]
	\arrow[from=2-1, to=2-3, "\cong"]
\end{tikzcd}
\end{equation*}
The left square arises from a square of ABC morphisms:
\[ \begin{tikzcd}[column sep = large]
{(\T,\I,F,A)} \arrow[d, "{[\Phi, \alpha, f]}"]                                                       && {(\T,\I,F,LA)} \arrow[ll, "{[\id_\T, \id_F, \mu_A]}"'] \arrow[d, "{\fk L(\fk m)}"] \\
{(\T',\I',F',A')} && {(\T',\I',F',L'A')} \arrow[ll, "{[\id_{\T'}, \id_{F'}, \mu'_{A'} ]}"']                        
\end{tikzcd} \]
Commutativity of this square boils down to the equation $\Phi(\mu'_{A'}) \circ \bb L_\Phi(A') \circ L(f) = f \circ \mu_A$, which follows from naturality of $\mu$ applied to $f$ and $\Phi \mu' \circ \bb L_\Phi = \mu \Phi$ applied to $A'$. The middle square commutes by \Cref{naturality of limit sheet phi prop}. Finally, the right square commutes as $\bb L(\fk m) = \App ( \fk L( \fk m))$.
\end{proof}
\end{thm*}

\section{Isomorphisms of K-theory groups}\label{applications section}

A morphism of spectral sequences which is an isomorphism on the second sheet will be an isomorphism on all later sheets and when the sequences converge, any compatible map of the target groups must be an isomorphism \cite[Comparison Theorem 5.2.12]{Weibel94}. An ABC morphism $\fk m \colon \fk M \to \fk M'$ which induces an isomorphism 
\[\bb L_n(\fk m) \colon \bb L_n(\fk M) \to \bb L_n(\fk M')\]
of derived functors for each $n \geq 0$ will therefore by \Cref{ABCthmintext} induce an isomorphism 
\[\bb L(\fk m) \colon \bb L(\fk M) \to \bb L(\fk M')\]
of localisations. Consider an ABC morphism $\fk m = [\Ind_\Omega, \alpha_\Omega, f]$ arising from an étale correspondence as in \Cref{motivating ABC morphism example}. If the families of open embeddings satisfy condition (P) then by \Cref{properly weak equivalence result} the localisations are the topological K-theories and $\bb L(\fk m) = \Ktop_*(\Omega;f)$ as in \Cref{functoriality of topological K theory}. Piecing together the above with the Baum--Connes conjecture we obtain a tool for proving isomorphisms of the K-theory groups of reduced crossed products:

\begin{cor}\label{main isomorphism corollary for baum connes}
Let $G$ and $H$ be second countable Hausdorff \'etale groupoids with coefficient C*-algebras $A \in \KK^G$ and $B \in \KK^H$, let $\Omega \colon G \to H$ be a second countable Hausdorff \'etale correspondence and let $f \in \KK^G(A,\Ind_\Omega B)$. Let $\cal E$ and $\cal F$ be families of second countable Hausdorff proper groupoids openly Morita embedded into $G$ and $H$ respectively which are compatible under $\Omega$, so that we get an ABC morphism 
\begin{equation*}
[\Ind_\Omega, \alpha_\Omega, f] \colon (\KK^G, \I_{\cal F_G},\K_*(G \ltimes -),A) \to (\KK^H, \I_{\cal F_H},\K_*(H \ltimes -),B).
\end{equation*}
Suppose that for each $n \geq 0$ the derived functor map
\begin{equation*}
\bb L_n(\Ind_\Omega,\alpha_\Omega,f) \colon \bb L_n \K_*(G \ltimes A) \to \bb L_n \K_*(H \ltimes B)
\end{equation*}
is an isomorphism. Then 
\[\Ktop_*(\Omega;f) \colon \Ktop_*(G;A) \to \Ktop_*(H;B)\]
is an isomorphism. If $(G,A)$ and $(H,B)$ satisfy the Baum--Connes conjecture we obtain an isomorphism $\K_*(G \ltimes_r A) \cong \K_*(H \ltimes_r B)$. If moreover the $\cs$-correspondence $\Omega \ltimes B \colon G \ltimes \Ind_\Omega B \to H \ltimes B$ descends to the reduced level, the resulting map $\Kred_*(\Omega;f) \colon \K_*(G \ltimes_r A) \to \K_*(H \ltimes_r B)$ equals this isomorphism (see \Cref{functoriality of topological K theory}).
\end{cor}

Application of this result depends on an understanding of the derived functors $\bb L_n$. Our understanding is especially good for ample groupoids with the singleton unit space family. The analysis in Section \ref{groupoid homology section} culminates in the description of $\bb L_n(\Ind_\Omega, \alpha_\Omega, f)$ in \Cref{induced maps in homology agree} as the induced map $\H_n(\Omega;\K_*(f))$ in homology from \cite[Theorem 3.5]{Miller24}. If the groupoids have torsion-free isotropy groups, the unit space family satisfies (P) by \Cref{checkable P}. Thus we may turn isomorphisms in homology into isomorphisms in K-theory: 
\begin{cor}
Let $G$ and $H$ be second countable Hausdorff ample groupoids with torsion-free isotropy groups and which satisfy the Baum--Connes conjecture with coefficients in $\cs$-algebras $A \in \KK^G$ and $B \in \KK^H$. Let $\Omega \colon G \to H$ be a second countable Hausdorff \'etale correspondence and let $f \in \KK^G(A,\Ind_\Omega B)$. Suppose that the induced map in homology 
\[\H_n(\Omega; \K_*(f))\colon \H_n(G; \K_*(A)) \to \H_n(H; \K_*(B))\]
is an isomorphism for each $n \geq 0$. Then there is an isomorphism $\K_*(G \ltimes_r A) \cong \K_*(H \ltimes_r B)$. 
\end{cor}
When $\Omega \colon G \to H$ is proper, taking the morphism $f_\Omega$ induced by the proper map $\orho \colon \Omega/H \to G^0$ gives us a version of this corollary without coefficients. Note that if $G$ and $H$ satisfy the Baum--Connes conjecture and $C^*(\Omega) \colon C^*(G) \to C^*(H)$ descends to a $\cs$-correspondence $C^*_r(\Omega) \colon C^*_r(G) \to C^*_r(H)$ at the reduced level, then the induced map in K-theory is always $\Kred_*(\Omega;f_\Omega) = \K_*(C^*_r(\Omega))$.
\begin{cor}\label{isomorphism of homology into k theory}
Let $G$ and $H$ be second countable Hausdorff ample groupoids with torsion-free isotropy groups and which satisfy the Baum--Connes conjecture and let $\Omega\colon G \to H$ be a proper second countable Hausdorff étale correspondence. Suppose that the induced map in homology $\H_{*}(\Omega)\colon \H_*(G) \to \H_*(H)$ is an isomorphism. Then there is an isomorphism $\K_*(C^*_r(G)) \cong \K_*(C^*_r(H))$. 
\end{cor}

Let $S$ be an inverse semigroup (with $0$) with idempotent semilattice $E$. There is a proper étale correspondence 
\[ \Omega_S \colon S^\times \to G_S \]
which induces an isomorphism in homology groups \cite[Example 3.10]{Miller24}. Here $S^\times$ is the discrete groupoid of non-zero elements of $S$, which has unit space $E^\times$, and $G_S$ is the universal groupoid. 

When $G_S$ is Hausdorff it is a groupoid model for the $\cs$-algebra $C^*_\lambda(S)$ generated by the left regular representation of $S$ by partial isometries on $\ell^2(S^\times)$ \cite[Theorem 5.5.18]{CELY17}. Hausdorffness of the universal groupoid $G_S$ is equivalent to the natural order on $S$ being a weak semilattice \cite[Theorem 5.17]{Steinberg10}. The isotropy group of $S^\times$ at an idempotent $e \in E^\times$ is the stabiliser subgroup $S_e = \{ s \in S \mid s^* s = s s^* = e \}$. Each stabiliser subgroup embeds in an isotropy group of $G_S$. Therefore if $G_S$ satisfies the Baum--Connes conjecture then so does $S^\times$ as a discrete groupoid.

Applying \Cref{isomorphism of homology into k theory} to the isomorphism $\H_*(\Omega_S)$ in homology we obtain a K-theory formula previously only achieved for inverse semigroups admitting an idempotent partial pure homomorphism into a group satisfying Baum--Connes \cite[Theorem 1.1]{Li21a}, at the price of adding a torsion-freeness hypothesis:
\begin{cor}\label{inverse semigroup corollary}
Let $S$ be a countable inverse semigroup with idempotent semilattice $E$ such that for each $e \in E^\times$ the stabiliser subgroup $S_e$ is torsion-free and the universal groupoid $G_S$ is Hausdorff and satisfies the Baum--Connes conjecture. Then
\[ \K_*(C^*_\lambda(S)) \cong \bigoplus_{\orb(e) \in S \backslash E^\times} \K_*(C^*_\lambda(S_e)). \]
\end{cor}

\begin{rmk}
In principle, this formula should have nothing to do with torsion-freeness of the stabiliser subgroups. In the presence of torsion, groupoid homology no longer describes the relevant derived functors, so we cannot apply \Cref{isomorphism of homology into k theory} which is itself built upon \Cref{groupoid homology section}. Instead, to apply our machinery to $\Omega_S \colon S^\times \to G_S$ we must work with more complex families of proper open subgroupoids of $S^\times$ and $G_S$. Following \Cref{inverse semigroup equivariant topological correspondence example}, the finite subgroups of stabiliser subgroups of $S$ determine families which are compatible under $\Omega_S \colon S^\times \to G_S$. \Cref{main isomorphism corollary for baum connes} then reduces the problem to showing that the associated ABC morphism $\fk m_S = [\Ind_{\Omega_S}, \alpha_{\Omega_S}, f_{\Omega_S}]$ induced by $\Omega_S$ induces an isomorphism $\bb L_n(\fk m_S)$ of derived functors for each $n \geq 0$. Consideration of all the finite subgroups of the stabiliser subgroups adds considerable complexity to this task which is therefore outside the scope of this article, but the torsion-free case serves as proof of concept for the approach.
\end{rmk}

\begin{example}[Right LCM monoids]\label{monoid example}
The left regular algebra $C^*_\lambda(P)$ of a countable right LCM left cancellative monoid $P$ is isomorphic to the left regular algebra $C^*_\lambda(I_{\ell}(P))$ of its left inverse hull $I_{\ell}(P)$. Assuming that the universal groupoid is Hausdorff, satisfies Baum--Connes and has torsion-free stabilisers, \Cref{inverse semigroup corollary} says that 
\[ \K_*(C^*_\lambda(P)) \cong \K_*(C^*_\lambda(P^*)), \]
where $P^*$ is the group of invertible elements of $P$. The main advancement from \cite[Corollary 1.4]{Li21a} is that $P$ need not embed in a group. Weakening the right LCM condition to the independence condition, there is a similar computation that incorporates the constructible ideals of $P$ (see \cite[Corollary 1.3]{Li21a}).
\end{example}

\begin{example}[The Toeplitz algebra of a singly aligned higher rank graph]\label{HRG example}
The Toeplitz algebra $\mathcal T C^*(\Lambda)$ of a countable singly aligned higher rank graph $\Lambda$ is isomorphic to the left regular algebra $C^*_\lambda(S_\Lambda)$ of the graph inverse semigroup $S_\Lambda$. The universal groupoid is Hausdorff, amenable, and has torsion-free stabilisers. The non-zero idempotents of $S_\Lambda$ may be identified with $\Lambda$ and their orbits under the action of $S_\Lambda$ correspond to the vertices $\Lambda^0$. The stabiliser subgroups are trivial, so \Cref{inverse semigroup corollary} recovers the K-theory computation \cite[Corollary 2.5.7]{Fletcher18}  
\begin{align*}
\K_0(\mathcal T C^*(\Lambda)) & \cong \bigoplus_{v \in \Lambda^0} \Z & \K_1(\mathcal T C^*(\Lambda)) & \cong 0.
\end{align*}
As higher rank graphs need not admit a faithful functor to a group (see \cite[Section 7]{PQR04}), the inverse semigroup $S_\Lambda$ may not admit an idempotent pure partial homomorphism to a group, so this computation is not covered by \cite[Theorem 1.1]{Li21a}.
\end{example}

%

\end{document}